\documentclass[12pt]{amsart}
\usepackage{amsthm}
\usepackage{amsmath, mathrsfs}
\usepackage{color}
\usepackage{amssymb}
\usepackage{hyperref}
\usepackage{enumerate}
\usepackage{latexsym,array}
\usepackage{amsfonts}
\usepackage{shadow}

\newtheorem{Pa}{Paper}[section]
\newtheorem{Tm}[Pa]{{\bf Theorem}}
\newtheorem{La}[Pa]{{\bf Lemma}}
\newtheorem{Dn}[Pa]{{\bf Definition}}

\newtheorem{Rk}[Pa]{{\bf Remark}}
\newtheorem{Pn}[Pa]{{\bf Proposition}}

\author[D. Alpay]{Daniel Alpay}
\address{Department of Mathematics\\
Ben–Gurion University of the Negev\\
Beer-Sheva 84105
Israel\\
dany@math.bgu.ac.il}
\author[F. Colombo]{Fabrizio Colombo}
\address{Politecnico di
Milano\\Dipartimento di Matematica\\Via E. Bonardi, 9\\20133 Milano,
Italy\\fabrizio.colombo@polimi.it}
\author[I. Sabadini]{Irene Sabadini}
\address{Politecnico di
Milano\\Dipartimento di Matematica\\Via E. Bonardi, 9\\20133 Milano,
Italy\\irene.sabadini@polimi.it}
\title[Schur functions in the slice hyperholomorphic setting]{
Schur functions and their realizations in the slice
hyperholomorphic setting} \oddsidemargin 0.2in \evensidemargin
0.2in \topmargin -0.5in \textwidth 15.5truecm \textheight 23truecm

\def\R{\mathbb R}

\def\s{\mathbf s}

\def\(s){\mathscr S(\R\times\R)}

\keywords{Schur functions, realization, reproducing
kernels, slice hyperholomorphic functions, $S$-resolvent
operators.}

\subjclass{MSC: 47B32, 47S10, 30G35}

\thanks{D. Alpay thanks the Earl Katz family for endowing the chair
which supported his research, and the Binational Science
Foundation Grant number 2010117.}
\begin{document}
\maketitle \tableofcontents
\parindent 0cm

\begin{abstract}
In this paper we start the study of
Schur analysis in the quaternionic setting using the
theory of slice hyperholomorphic
functions. The novelty of our approach is that slice hyperholomorphic
functions allows to write realizations in terms of a suitable resolvent,
the so called S-resolvent operator and to extend several results that hold
in the complex case to the quaternionic case. We discuss reproducing kernels,
positive definite functions in this setting and we show how they can be
obtained in our setting using the extension operator and the slice regular
product. We define Schur multipliers, and find their co-isometric realization
in terms of the associated de Branges-Rovnyak space.
\end{abstract}

\section{Introduction}
\setcounter{equation}{0}
In this paper we develop Schur analysis, and in particular the
Schur algorithm, and a theory of linear systems when the complex
numbers are replaced by the skew-field of quaternions.  An
important tool is the theory of slice hyperholomorphic
functions. So there is a combination of a non-commutative setting (since
the quaternions lack the commutativity property) and of analyticity (via the
slice hyperholomorphic functions).
Since the paper is aimed at two different audiences,
namely researchers from Clifford analysis and researchers from
operator theory and classical linear system theory,
we will survey the basic definitions from both fields needed in the paper.\\

We denote by $\mathcal S$ the set of functions analytic and
contractive in the open unit disk $\mathbb D\subset\mathbb{C}$. Such functions
bear various names, and we will call them {\sl Schur functions}
in the present paper. Let $s\in\mathcal S$, and assume that
$|s(0)|<1$ (and therefore, by the maximum modulus principle, $s$
is not equal to a unitary constant, but takes strictly
contractive values in $\mathbb D$). Then, the function
\begin{equation}
\label{voltaire}
s^{(1)}(z)=\begin{cases}\dfrac{1}{z}\dfrac{s(z)-s(0)}{1-s(z)\overline{s(0)}},
\quad z\not =0,\\
\dfrac{s^{\prime}(0)}{1-|s(0)|^2},\quad \hspace{0.6cm}z=0,
\end{cases}
\end{equation}
also belongs to $\mathcal S$. More generally, the recursion
\begin{equation}
\label{recur}
\begin{split}
s^{(0)}&=s(z)\\
s^{(n+1)}(z)&=\begin{cases}\dfrac{1}{z}\dfrac{s^{(n)}(z)-s^{(n)}(0)}
{1-s^{(n)}(z)\overline{s^{(n)}(0)}},
\quad z\not =0,\\
\dfrac{(s^{(n)})^{\prime}(0)}{1-|s^{(n)}(0)|^2},\quad
\hspace{1.cm}z=0,
\end{cases}
\end{split}
\end{equation}
defines a sequence, finite or infinite, of Schur functions
$s^{(0)},s^{(1)},\ldots$. The sequence is infinite if
\[
s^{(n)}(0)\in\mathbb D,\quad n=0,1,\ldots,
\]
while it stops at rank $n$ if $|s^{(n)}(0)|=1$. The construction
of this sequence is the celebrated {\sl Schur algorithm},
developed by I. Schur in 1917. See \cite{schur, schur2,goh1}. The
numbers $\rho_n=s^{(n)}(0)$ are called the Schur coefficients
associated to $s$, and the sequence (finite or infinite) of Schur
coefficients uniquely determines $s$.\\

Let us already make the following remark at this point. In the
Schur algorithm one makes use of Schwarz's lemma and of the
elementary fact that if $u$ and $v$ are in the open unit disk, so
is
\begin{equation}
\label{uv} \frac{u-v}{1-u\overline{v}}.
\end{equation}
Let us denote by $\mathbb B$ the open unit ball of the
quaternions. If one replaces in \eqref{uv} $u$ and $v$ by
quaternions of norm strictly less than $1$, then both
\[
(u-v)(1-u\overline{v})^{-1}\quad {\rm and}\quad
(1-u\overline{v})^{-1}(u-v)
\]
are still in $\mathbb B$. But, as we shall see, these transformations will not
keep the property of being slice hyperholomorphic, and a different
approach will be needed. Note that the Schwarz lemma will not
help to develop a Schur algorithm here because of the lack of
commutativity. If $s$ is slice hyperholomorphic into $\mathbb B$
(and in particular $|s(0)|<1$), the functions
\[
(1-s(p)\overline{s(0)})^{-1}(s(p)-s(0)),
\]
or
\[
(s(p)-s(0))(1-s(p)\overline{s(0)})^{-1}
\]
are indeed contractive in $\mathbb B$, but they will not be
slice hyperholomorphic. So there is no direct counterpart of
\eqref{voltaire}. One needs to use the notion of Schur multipliers
and the slice regular product.\\

{\sl Schur analysis} originates with the works of Schur, Herglotz, and others
(see \cite{hspnw} for reprints of the original works), and can be seen as
a collection of questions pertaining to Schur functions and their various
generalizations to other settings. Among the problems we mention in
particular:\\
$(1)$ Classical interpolation problems such as
Carath\'eodory-Fej\'er and Nevanlinna-Pick problems, and their
matrix-valued versions; see for instance \cite{dfk,Dym_CBMS}. The
case of boundary interpolation is of special importance. See for
instance
\cite{MR2775147} for related recent results.\\
$(2)$ Realization of Schur functions in the form
$s(z)=D+zC(I-zA)^{-1}B$, where the operator matrix
\[
\begin{pmatrix}A&B\\ C&D\end{pmatrix}
\]
is subject to various metric constraints, namely coisometric, isometric and
unitary. See \cite{dbr1, dbr2, adrs}.\\
$(3)$ Schur functions are closely related to the theory of linear
systems. The term {\sl linear system} encompasses a wide range of
situations. Here we have in mind input output of the form
\[
y_n=\sum_{m=0}^ns_m u_{n-m},\quad n=0,1,\ldots,
\]
where $s_0,s_1,\ldots$ is a sequence of matrices of $\mathbb
C^{M\times N}$ (the impulse response), $(u_n)_{n\in\mathbb N_0}$
is a sequence of vectors of $\mathbb C^N$ (the input sequence),
and $(y_n)_{n\in\mathbb N_0}$ is a sequence of vectors of
$\mathbb C^M$ (the output sequence). The function
$s(z)=\sum_{n=0}^\infty s_nz^n$ is a Schur function if and only if
the $\ell_2$ norm of the output sequence is always less or equal
to the $\ell_2$ norm of the input sequence. In other words, Schur
functions are the transfer functions of time-invariant dissipative
linear systems.\\
$(4)$ Yet another direction of research is related to inverse
scattering; see
\cite{MR88g:94024,MR2003a:47139,ad1, ad2}.\\
$(5)$ Last but not least we mention the connection with fast algorithms; see
\cite{SK, SK1}.

These various questions make sense in more general settings, of
which we mention in particular the several complex variables
case, the indefinite case, the time-varying case, the
non-commutative case, the case of compact Riemann surfaces, and
the case of several complex variables, to name a few. References
are given in the last section of the paper.  In the present work
we present a counterpart of Schur functions and of the Schur
algorithm in the quaternionic setting.\\

The main tool that we use to extend Schur analysis to the
quaternionic setting is the theory of slice hyperholomorphic, or slice regular, functions.
 Some references for this theory of functions are \cite{MR2353257, MR2555912, MR2742644}. For the generalization to functions with values in a Clifford algebra, called again slice hyperholomorphic or slice monogenic functions, we refer the reader to \cite{MR2520116, MR2684426, MR2673423}, and to \cite{MR2737796} for functions with values in a real alternative algebra. Finally, we mention that it exists a non constant coefficients differential operator whose kernel contains slice hyperholomorphic
 functions defined on suitable domains, \cite{GLOBAL}.
The theory of slice regular functions allows to define the
quaternionic functional calculus and its associated $S$-resolvent
operator. The importance of the $S$-resolvent operator is in the
definition of the quaternionic version of the operator
$(I-zA)^{-1}$ that appears in the realization function
$s(z)=D+zC(I-zA)^{-1}B$. It turns out that when $A$ is a
quaternionic matrix and $p$ is a quaternion then $(I-pA)^{-1}$
has to be replaced by $( I  -\bar p A)(|p|^2A^2-2{\rm Re}(p) A+
I  )^{-1}$ which is equal to $p^{-1}S^{-1}_R(p^{-1},A)$ where
$S^{-1}_R(p^{-1},A)$ is the right $S$-resolvent operator
associated to the quaternionic matrix $A$. For some results on
the quaternionic functional calculus we refer the reader to
\cite{MR2735309, MR2496568, MR2661152, MR2803786}. Slice
monogenic functions admit a functional calculus for $n$-tuples of
operators and for this theory we mention \cite{MR2402108,
MR2720712, SCFUNCTIONAL}. The book \cite{MR2752913} collects some
of the main results on the theory of slice hyperholomorphic
functions and the related functional calculi.
\\

Finally we mention the paper \cite{MR2124899,MR2240272,asv-cras},
where Schur multipliers were introduced and studied in the quaternionic
setting using the
Cauchy-Kovalesvkaya product and series of Fueter polynomials, and
the papers \cite{MR733953,MR2205693,MR2131920}, which treat
various aspects of a theory of linear systems in the quaternionic
setting. Our approach is quite different from the methods used
there. \\

The paper consists of eight sections besides the introduction, and
its outline is as follows: in Section \ref{sec:1} we review the
main aspects of the theory of slice hyperholomorphic functions
and the $S$-resolvent operators. In Section \ref{sec2} we study
the counterpart of state space equations in the slice
hyperholomorphic setting. This leads us naturally to the notion
of rational function, defined and studied in Section
\ref{sec:real}. In Section \ref{sec:positive} we study certain
positive definite functions. This paves the way to the definition
of Schur multipliers and to a version of  Schwarz' lemma in
Section \ref{sec:schur_mul}. In Section \ref{sec:co} we exhibit a
coisometric realization of a Schur multiplier. Section
\ref{schur_h} is devoted to the Schur algorithm in the present
setting. In the last section we discuss briefly future directions
of research.

\section{Slice hyperholomorphic functions and the $S$-resolvent operators}
\setcounter{equation}{0}
\label{sec:1}

In this section we introduce the preliminary results related to the theory
of slice hyperholomorphic functions
and the quaternionic $S$-resolvent operators.
In the next section we will use these tools to define in a suitable way the
quaternionic slice hyperholomorphic transfer function associated to quaternionic linear systems.\\

By $\mathbb{H}$ we denote the algebra of real quaternions
$p=x_0+ix_1+jx_2+kx_3$ which can also we written as $p={\rm
Re}(p)+{\rm Im}(p)$ where $x_0={\rm Re}(p)$ and
$ix_1+jx_2+kx_3={\rm Im}(p)$. By $\mathbb{S}$ we indicate the set
of unit purely imaginary quaternions, i.e.
$$
\mathbb{S}=\{p=x_1i+x_2j+x_3k\ :\ x_1^2+x_2^2+x_3^2=1\}.
$$
In the literature there are several notions of quaternion valued
hyperholomorphic functions. In this paper we will consider a
notion which includes power series in the quaternionic variable,
the so-called slice regular or slice hyperholomorphic functions. The main reference for the material in this section is
the book \cite{MR2752913}.
\begin{Dn}
Let $\Omega\subseteq\mathbb{H}$ be an open set and let
$f:\ \Omega\to\mathbb{H}$ be a real differentiable function. Let
$I\in\mathbb{S}$ and let $f_I$ be the restriction of $f$ to the
complex plane $\mathbb{C}_I := \mathbb{R}+I\mathbb{R}$ passing through $1$
and $I$ and denote by $x+Iy$ an element on $\mathbb{C}_I$.
\begin{enumerate}
\item
 We say that $f$ is a left slice regular function
(or slice regular or slice hyperholomorphic)  if, for every
$I\in\mathbb{S}$, we have:
$$
\frac{1}{2}\left(\frac{\partial }{\partial x}+I\frac{\partial
}{\partial y}\right)f_I(x+Iy)=0.
$$
\item
We say that $f$ is right slice regular function (or right slice hyperholomorphic) if,
for every
$I\in\mathbb{S}$, we have
$$
\frac{1}{2}\left(\frac{\partial }{\partial x}f_I(x+Iy)+\frac{\partial
}{\partial y}f_I(x+Iy)I\right)=0.
$$
\end{enumerate}
\end{Dn}
\begin{Dn}
Given $p\in\mathbb{H}$, if $q$ is not real we can define
$I_p:={\rm Im}(p)/|{\rm Im}(p)|$, so $p={\rm Re}(p)+I_p |{\rm
Im}(p)|$. We denote by $[p]$ the set of all elements of the form
${\rm Re}(p)+J |{\rm Im}(p)|$ when $J$ varies in $\mathbb{S}$. We
say that $[p]$ is the 2-sphere defined by $q$.
\end{Dn}

\begin{Dn}
Let $\Omega$ be a domain in $\mathbb{H}$. We say that $\Omega$ is
a \textnormal{slice domain} (s-domain for short) if $\Omega \cap
\mathbb{R}$ is non empty and if $\Omega\cap \mathbb{C}_I$ is a
domain in $\mathbb{C}_I$ for all $I \in \mathbb{S}$. We say that
$\Omega$ is \textnormal{axially symmetric} if, for all $p\in
\Omega$, the 2-sphere $[p]$ is contained in $\Omega$.
\end{Dn}
In the sequel we will work mainly on the unit sphere in $\mathbb{H}$
with center at the origin, which is trivially an axially symmetric s-domain.
\begin{Tm}[Representation Formula]\label{formula}
Let $\Omega\subseteq \mathbb{H}$ be an axially symmetric s-domain.
Let
$f$ be a left slice regular function on  $\Omega\subseteq  \mathbb{H}$.
 Then the following equality holds for all $p=x+I_p y \in \Omega$:
\begin{equation}\label{strutturaquat}
f(p)=f(x+I_p y) =\frac{1}{2}\Big[   f(z)+f(\overline{z})\Big]
+\frac{1}{2}I_pI\Big[ f(\overline{z})-f(z)\Big],
\end{equation}
where $z:=x+Iy$, $\overline{z}:=x-Iy\in\mathbb{H}$. Let $f$ be a
right slice regular function on  $\Omega\subseteq  \mathbb{H}$.
Then the following equality holds for all $p=x+I_p y \in \Omega$:
\begin{equation}\label{distribution}
f(x+I_p y) =\frac{1}{2}\Big[   f(z)+f(\overline{z})\Big]
+\frac{1}{2}\Big[ f(\overline{z})-f(z)\Big]II_p.
\end{equation}
\end{Tm}
The Representation Formula allows to extend any function
$f: \ \Omega\subseteq\mathbb{C}_I\to\mathbb{H}$ defined on an axially
symmetric open set  $\Omega$ intersecting the real axis and in the kernel of the
Cauchy-Riemann operator to
a function $f: \ \widetilde{\Omega}\subseteq\mathbb{H}\to\mathbb{H}$ slice
regular where $\widetilde{\Omega}$ is the smallest axially symmetric
open set in $\mathbb{H}$ containing $\Omega$ by means of the {\em extension operator}
\begin{equation}\label{ext}
{\rm ext}(f)(p):= \frac{1}{2}\Big[   f(z)+f(\overline{z})\Big]
+\frac{1}{2}I_pI\Big[ f(\overline{z})-f(z)\Big],\quad z,\bar
z\in\mathbb{C}_I.
\end{equation}
Slice regular functions satisfy an identity principle,
specifically, two of them coincide in a domain $\Omega$
intersecting the real axis if they coincide on a subset of
$\Omega\cap\mathbb{R}$ having an accumulation point. As a consequence we have:
\begin{Pn}\label{realanalytic}
Any  quaternion valued real analytic function defined in
$(a,b)\subseteq\mathbb{R}$ extends uniquely to a slice regular function
defined in a suitable open set containing $(a,b)$.
\end{Pn}
\begin{proof}
Let $f:\ (a,b)\subseteq\mathbb{R}\to\mathbb{H}$ be a real analytic function.
Then it can be uniquely extended to a suitable domain $D$ containing $(a,b)$, by
extending its real components, to a function
$\tilde f:\ D\subseteq\mathbb{C}\to\mathbb{H}$ which is in the kernel of
the Cauchy-Riemann operator. By using the extension operator (\ref{ext}) $\tilde f$
extends to a slice hyperholomorphic function defined on the smallest axially
symmetric open set $\tilde\Omega$ containing $D$. The uniqueness follows from
the identity principle.
\end{proof}
Given two slice regular functions their product is not, in
general, slice regular. It is possible to introduce a suitable
product denoted by $\star$, see \cite[Definition 4.3.5, p.
125]{MR2752913} (note that in the literature this product is
denoted by $*$ but here we use this symbol in order to avoid
confusion with adjoint of operators). Here we describe the
$\star$-product which gives a left slice regular product. It is
possible to define an analog product to multiply two right slice
regular functions and obtain a function with the same regularity.
Let $\Omega\subseteq\mathbb{H}$ be an axially symmetric s-domain
and let
 $f,g:\ \Omega\to\mathbb{H}$ be slice regular functions. Fix $I,J\in\mathbb{S}$, with
 $I\perp J$. The Splitting Lemma guarantees the existence of four holomorphic functions
 $F,G,H,K: \ \Omega\cap\mathbb{C}_I\to \mathbb{C}_I$ such
 that for all $z=x+Iy\in   \Omega\cap\mathbb{C}_I$
 $$
 f_I(z)=F(z)+G(z)J, \qquad g_I(z)=H(z)+K(z)J.
 $$
Define the function $f_I\star g_I:\ \Omega\cap\mathbb{C}_I\to \mathbb{H}$ as
\begin{equation}
f_I\star g_I(z)=[F(z)H(z)-G(z)\overline{K(\bar z)}]
+[F(z)K(z)+G(z)\overline{H(\bar z)}]J.
\end{equation}
Then  $f_I\star g_I(z)$ is obviously a holomorphic map and hence
we can  give the following definition.
 \begin{Dn}
Let $\Omega\subseteq\mathbb{H}$ be an axially symmetric s-domain and let
 $f,g:\ \Omega\to\mathbb{H}$  be slice regular. The function
$$f\star g(p)={\rm ext}(f_I\star g_I)(p)$$
is called the slice regular product of $f$
and $g$.
 \end{Dn}
 \begin{Rk}{\rm
It is immediate to verify that the $\star$-product
is associative, distributive but, in general, it is not commutative.
}
\end{Rk}

When slice regular functions are defined on a ball
$\mathbb{B}_R\subseteq\mathbb{H}$ with center at the origin and
radius $R$ then they admit power series expansions. Suppose that
on $\mathbb{B}_R$ we have $f(p)=\sum_{n=0}^\infty p^n a_n$ and
$g(p)=\sum_{n=0}^\infty p^n b_n$, then
  the $\star$-product is given by
\begin{equation}
\label{eq:Cauchy_product}
(f\star g)(p)=\sum_{n=0}^\infty p^n c_n,
\qquad c_n=\sum_{r=0}^n a_rb_{n-r}.
\end{equation}
In other words, the coefficient sequence $(c_n)_{n\in\mathbb
N_0}$ is the Cauchy product, or the convolution product, of the
sequences of coefficients associated to $f$ and $g$. See for
instance \cite[(2) p. 199]{MR51:583} for the Cauchy product
in a non -commutative setting.
\begin{Rk}{\rm
From Proposition \ref{realanalytic}, it immediately follows that
$f\star g$ is uniquely determined by $f(x)g(x)$, $x\in
\mathbb{B}_R\cap\mathbb{R}$. }
\end{Rk}
We observe that another product in Clifford analysis (see
\cite{bds, MR2089988}), namely the Cauchy-Kovaleskaya product
which allows to obtain Cauchy-Fueter regular functions or
monogenic functions, can be seen as a (different) convolution.
More generally, pointwise product is often best replaced by
convolution of underlying coefficient sequences. See \cite{bds}.
We mention as an example the Wick product in white noise space
analysis. See \cite{MR1408433}, and see \cite{al_acap,alp}
for applications to linear system theory.\\

\begin{Rk}\label{gcomf}  Let $f(p)=\sum_n p^n a_n$ and $g(p)=\sum_n p^n b_n$.
If $f$ has real coefficients then $f\star g=g\star f$.
\end{Rk}
Given a slice regular function $f$ it is possible to construct its slice regular
reciprocal, denoted by $f^{-\star}$. We do not provide the
general construction, which can be found in \cite{MR2752913}
since it is sufficient to construct the inverse of a polynomial
or a power series with center at the origin. Given $f(q)$ as
above, let us introduce the notation
$$
f^c(p)=\sum_{n=0}^\infty p^n \bar a_n,\qquad  f^s(p)=(f^c\star f)(p
)=\sum_{n=0}^\infty p^nc_n,\quad
c_n=\sum_{r=0}^n a_r\bar a_{n-r}.
$$
Note that the series $f^s$ has real coefficients. The slice regular reciprocal
is then defined as
$$
f^{-\star}:=(f^s)^{-1}f^c.
$$
 In an analogous way on can define
the right slice regular reciprocal $ f^{-\star}:=f^c(f^s)^{-1}, $
of a right regular function  $f(p)=\sum_n a_n p^n $.

We do not introduce a specific symbol in order to distinguish the
left and the right slice regular reciprocal, since the context
will clarify the case we are considering. There is however a
remarkable case of reciprocal that deserves further explanations.
Consider the function
$$
S(r,p):=p-q,  \ \ \ \ \ p, q\in \mathbb{H}.
$$
Its slice regular reciprocal can be constructed in four possible
ways: we can construct a reciprocal which is left (resp. right)
regular with respect to the variable $p$ or left (resp. right)
regular with respect to the variable $q$. Accordingly to these
possibilities we obtain the function (see \cite{MR2752913})
$$
 S_L^{-1}(p,q)=-(q^2-2q {\rm Re}(p)+|p|^2)^{-1}(q-\overline{p})
$$
which corresponds to the left slice regular reciprocal  in the
variable $q$. The function $S_L^{-1}(p,q)$  is left regular in the
variable $q$ by construction and it turns out to be right regular
with respect to $p$. While the other possibility gives
$$
{S}_R^{-1}(p,q):=-(q-\bar p)(q^2-2{\rm Re}(p)q+|p|^2)^{-1},
$$
which is right slice regular in $q$ and left slice regular with
respect to $p$.
\begin{Rk}
\label{-star} {\rm When no confusion will arise we will write
instead of $(p-q)^{-\star}$ its explicit expression using
${S}_L^{-1}(p,q)$ or ${S}_R^{-1}(p,q)$, according to the left or
right slice regularity required.}
\end{Rk}

It is possible to show that both kernels can be written in two different ways
$$
 S_L^{-1}(p,q)=-(q^2-2q
{\rm Re}(p)+|p|^2)^{-1}(q-\overline{p})=(p-\bar q)(r^2-2{\rm
Re}(q)p+|q|^2)^{-1},
$$
and
$$
S_R^{-1}(p,q) =-(q-\bar p)(q^2- 2 {\rm Re}(p)q+|p|^2)^{-1}= (p^2-2
{\rm Re} (q)r+|q|^2)^{-1}(r-\bar q),
$$
thus we have
$$
S_R^{-1}(p,q) = - S_L^{-1}(q,p) .
$$

Given two polynomials $P(p)$, $Q(p)$ it is possible to define two
possible "quotients" as $P^{-\star}\star Q$ or  $Q\star
P^{-\star}$ and due to the noncommutativity of quaternions, these
quotients do not coincide. Thus it is necessary to make a choice
between the left and the right quotient of $P$ and $Q$ and this
will be one of the various definitions of {\sl rational
functions}. A number of equivalent characterization of rational
functions are
given in Section \ref{sec:real}.\\

The realization of Schur functions in the form
$s(z)=D+zC(I-zA)^{-1}B$, where one wishes to replace the complex
variable $z$ by a quaternionic variable, and where the operator
matrix
\[
\begin{pmatrix}A&B\\ C&D\end{pmatrix}
\]
is now  a quaternionic operator matrix, requires a new concept to
replace the classical resolvent operator $(I-zA)^{-1}$. This new
object is the so called $S$-resolvent operator associated to the
quaternionic functional calculus. Let $V$ be a two sided
quaternionic Banach space, we denote by $\mathcal{B}(V)$ the
quaternionic Banach spaces of all bounded linear (left or right)
operators  endowed with the  natural norm. A remarkable fact is
that the following theorem  holds in the case in which the
components of the quaternionic operator $A$ do not commute.
\begin{Tm}\label{Ssinistro}
Let $A\in\mathcal{B}(V)$.
Then,  for $\|A\|< |r|$,  we have
\begin{equation}\label{ciaoleft}
\sum_{n= 0}^\infty A^n r^{-1-n}=-
(A^2-2{\rm Re}(r) A+|r|^2 I  )^{-1}(A-\overline{r} I  ),
\end{equation}
\begin{equation}\label{SresolvR}
\sum_{n= 0}^\infty  r^{-1-n}A^n=-
(A- I  \overline{r})(A^2-2{\rm Re}(r) A+|r|^2 I  )^{-1},
\end{equation}
where $ I  $ denotes the identity operator.
\end{Tm}

The notion of $S$-spectrum of a quaternionic operator $A$ is
suggested by the definition of $S$-resolvent operator that is the
kernel for the quaternionic functional calculus; see
\eqref{ciaoleft} and \eqref{SresolvR}.

\begin{Dn}[The $S$-spectrum and the $S$-resolvent sets of  quaternionic operators]
\label{defspscandres}
Let $T\in\mathcal{B}(V)$.
We define the $S$-spectrum $\sigma_S(A)$ of $T$  as:
$$
\sigma_S(A)=\{ r\in \mathbb{H}\ \ :\ \ A^2-2 \ {\rm Re}(r)A+|r|^2 I  \ \ \
{\it is\ not\  invertible}\}.
$$
The $S$-resolvent set $\rho_S(A)$ is defined by
$$
\rho_S(A)=\mathbb{H}\setminus\sigma_S(A).
$$
\end{Dn}

\begin{Rk}\label{comstrucs}{\rm
We have proved that
the $S$-spectrum
 $\sigma_S (A)$  is a compact nonempty set and if $p\in   \mathbb{H}$ belongs to $\sigma_S(A)$, then all the elements of the sphere $[p]$
 are contained in  $\sigma_S(A)$.
}
\end{Rk}

\begin{Dn}[The $S$-resolvent operator]
Let $V$ be a bilateral quaternionic Banach space, $A\in\mathcal{B}(V)$ and $r\in \rho_S(A)$.
We define the left $S$-resolvent operator as
\begin{equation}\label{Sresolvoperatordd}
S_L^{-1}(r,A):=-(A^2-2{\rm Re}(r) A+|r|^2 I  )^{-1}(A-\overline{r} I  ).
\end{equation}
We define the right $S$-resolvent operator as
\begin{equation}\label{SresolvoperatorRdd}
S_R^{-1}(r,A):=-(A- I  \overline{r})(A^2-2{\rm Re}(r) A+|r|^2 I  )^{-1}.
\end{equation}
\end{Dn}

\begin{Tm}[The $S$-resolvent equation]
Let $V$ be a bilateral quaternionic Banach space, $A\in\mathcal{B}(V)$ and  $r\in \rho_S(A)$.
Then the left $S$-resolvent
operator defined in (\ref{Sresolvoperatordd}) satisfies the equation
\begin{equation}\label{resSLEFT}
S_L^{-1}(r,A)r-AS_L^{-1}(r,A)= I  .
\end{equation}
and the right $S$-resolvent
operator defined in (\ref{SresolvoperatorRdd}) satisfies the equation
\begin{equation}\label{resSRIGHT}
sS_R^{-1}(s,A)-S_R^{-1}(s,A)A= I  .
\end{equation}
Moreover, $S_L^{-1}(r,\cdot): \rho_S(A)\to \mathcal{B}(V)$ is right slice regular in the variable $r$ and
$S_R^{-1}(r,\cdot): \rho_S(A)\to \mathcal{B}(V)$ is left slice regular in the variable $r$.
\end{Tm}

The notation introduced in Remark \ref{-star} is not necessarily valid when we replacing operators in place of the quaternionic variables.
We prove that it is still the case in the
 following proposition which is crucial
  and can be proved by suitably modifying the proof of Theorem \ref{Ssinistro}.
\begin{Pn}
Let $\mathcal H$ be a two sided quaternionic Hilbert space and let
$A$ be a bounded right linear quaternionic operator from $\mathcal H$ into
itself. Then, for $|p| \,\|A\|< 1$ we have
\begin{equation}
\label{eq:oberkampf_ligne_5}
\sum_{n=0}^\infty p^n A^n =p^{-1}S^{-1}_R(p^{-1},A)=( I  -\bar p A)(|p|^2A^2-2 {\rm Re}(p) A+ I  )^{-1},
\end{equation}
and
\begin{equation}\label{eq:oberkampf_ligne_55}
(I-pA)^{-\star}=\sum_{n=0}^\infty p^n A^n.
\end{equation}
\end{Pn}

\begin{proof} From (\ref{SresolvR}), it is immediate that
\[
\sum_{n=0}^\infty p^n A^n=p^{-1}S^{-1}_R(p^{-1},A)
\]
and by writing explicitly $S^{-1}_R(p^{-1},A)$ we obtain
$$
p^{-1}S^{-1}_R(p^{-1},A)
=-p^{-1}(A-\frac{p}{|p|^2} I  )(A^2-2\frac{{\rm
Re}(p)}{|p|^2}A+\frac{1}{|p|^2} I  )^{-1}
$$
and with some computations we get the second equality in \eqref{eq:oberkampf_ligne_5}.
To prove (\ref{eq:oberkampf_ligne_55}) we consider the function $p^{-1}-q$ and its slice regular reciprocal with respect to $q$ which is obtained  using the formula $f^{-\star}=f^c(f^s)^{-1}$, see Section 2. We obtain  $(p^{-1}-q)^{-\star}=S^{-1}_R(p^{-1},q)$ and, by the functional calculus, we can substitute $q$ by a quaternionic operator $A$ (it is crucial to observe that the components of $A$ do not necessarily commute) and we obtain:
$$
(p^{-1} I  -A)^{-\star}=-S_R^{-1}(p^{-1},A)=-(A-\bar p^{-1} I  )(A^2-2{\rm
Re}(p^{-1})A+|p^{-1}|^2 I  )^{-1}
$$
since
\[
p^{-1}(p^{-1} I  -A)^{-\star}=[p(p^{-1} I  -A)]^{-\star}=( I  -pA)^{-\star}.
\]
By using (\ref{eq:oberkampf_ligne_5}) we obtain (\ref{eq:oberkampf_ligne_55}).
\end{proof}

\section{State space equations and realization}
\setcounter{equation}{0}
\label{sec2}
We now show that if we consider   the  quaternionic linear system
\begin{equation}\label{linquatsy}
\left\{\begin{array}{l}
x_{n+1}=x_nA+u_nB,\ \ \ n=0,1,...\\
y_n=x_nC+u_nD,\\
\end{array}
\right.
\end{equation}
where  $A$, $B$, $C$, $D$ are given quaternion matrices, then the
quaternionic ''transfer function"  cannot be defined by simply
taking $\mathcal Z$-transform as in the complex case. The
following considerations show the problem that arises: Given a
sequence of quaternions (or of quaternionic matrices of common
size) $U=\{u_n\}_{n\in \mathbb{N}}$, we define the quaternionic
$\mathcal{Z}$-transform as
$$
\mathcal{Z}(U):=\mathcal{U}(p):=\sum_{n=0}^{\infty}p^nu_n,
$$
so the $\mathcal{Z}$-transform is right linear
$$
\mathcal{Z}(UA)=\mathcal{Z}(U)A.
$$
Since $\mathcal{Z}(U)$ is a power series centered at the origin it
is slice hyperholomorphic. The main properties of the classical
$\mathcal{Z}(U)$-transform still hold in the quaternionic setting. In
particular, if we set
$$
\tau_{-1}U:=(u_1,u_2,\ldots)
$$
we have
$$
\mathcal{Z}(\tau_{-1}U)=p^{-1}\mathcal{Z}(U) \ \ {\rm if} \ \ \ u_0=0.
$$
If we apply the $\mathcal{Z}$-transform to system (\ref{linquatsy}) we get
\begin{equation}\label{linquatsy2}
\left\{\begin{array}{l}
p^{-1}\mathcal{X}(p)= \mathcal{X}(p)A+\mathcal{U}(p)B,\ \ \ n=0,1,...\\
\mathcal{Y}(p)=\mathcal{X}(p)C+\mathcal{U}(p)D.\\
\end{array}
\right.
\end{equation}
The natural definition  of the transfer function of the system is
$$
\mathcal{H}(p):=(\mathcal{U}(p))^{-1}\mathcal{Y}(p)
$$
so we have to solve the quaternionic equation
$$
p^{-1}\mathcal{X}(p)- \mathcal{X}(p)A=B\mathcal{U}(p)
$$
which has the solution
$$
\mathcal{X}(p)=\sum_{n=0}^{\infty}p^nB\mathcal{U}(p)A^{n+1}.
$$
Because of the term $B$, this expression need not be slice
hyperholomorphic. By replacing this solution in the second
equation in (\ref{linquatsy2}), we obtain
$$
\mathcal{Y}(p)=\sum_{n=0}^{\infty}p^nB\mathcal{U}(p)A^{n+1}C+\mathcal{U}(p)D.
$$
It turns out that $\mathcal{H}(p)$ depends on $\mathcal{U}(p)$, in fact
$$
\mathcal{H}(p):=(\mathcal{U}(p))^{-1}\big(\sum_{n=0}^{\infty}p^nB\mathcal{U}(p)A^{n+1}
C+\mathcal{U}(p)D\Big)
$$
$$
=
(\mathcal{U}(p))^{-1}\sum_{n=0}^{\infty}p^nB\mathcal{U}(p)A^{n+1}C+D,
$$
and need not be slice hyperholomorphic. We now follow  an approach
based on slice hyperholomorphic to overcome the above difficulty
and to give a good
 definition of transfer function.

%Here we consider the complex case but we use the
% the extension operator (\ref{ext}) and
%the $\star$-product in order to obtain a slice regular function.

\begin{Tm}
Let $A$, $B$, $C$, $D$ be given  matrices of appropriate size
with quaternionic entries. Suppose that $\{u_n\}_{n\in
\mathbb{N}_0}$ is a given sequence of vectors with quaternionic
entries, and of appropriate size. Consider the system
\begin{equation}\label{linquatsydffrg}
\left\{\begin{array}{l}
x_{n+1}=Ax_n+B u_n,\ \ \ n=0,1,...\\
y_n=Cx_n+Du_n,\\
\end{array}
\right.
\end{equation}
and define its slice hyperholomorphic transfer function matrix-valued as
$$
\mathcal{H}(p):=  \mathcal{Y}(p)\star(\mathcal{U}(p))^{-\star}
$$
where $\mathcal{Y}(p)$ and $\mathcal{U}(p)$ are the slice hyperholomorphic extensions of the $Z$-transforms of $y_n$ and of $u_n$, respectively.
Then we have
\begin{equation}\label{tranfslice}
\mathcal{H}(p)= D+pC\star (I-pA)^{-*}B.
\end{equation}
Moreover, in terms of the right $S$-resolvent operator (\ref{tranfslice}) can be written as
\begin{equation}\label{reso}
\mathcal{H}(p)=C\ \star S^{-1}_R(p^{-1},A) B + D.
\end{equation}
\end{Tm}
\begin{proof}
Let us first consider system (\ref{linquatsydffrg}) in the complex plane $\mathbb{C}_I$, for a
fixed $I\in \mathbb{S}$, where
$A$, $B$, $C$, $D$,  now denoted by $a$, $b$, $c$, $d$, are given element in $\mathbb{C}_I$.
Suppose that $\{u_n\}_{n\in \mathbb{N}_0}$ is
a given sequence of complex numbers in $\mathbb{C}_I$:
\begin{equation}\label{linq}
\left\{\begin{array}{l}
x_{n+1}=ax_n+b u_n,\ \ \ n=0,1,...\\
y_n=cx_n+du_n.\\
\end{array}
\right.
\end{equation}
By setting
$$
z:=u+Iv\in \mathbb{C}_I
$$
and by taking the classical $\mathcal{Z}$-transform we get
\begin{equation}\label{linqu3}
\left\{\begin{array}{l}
\mathcal{X}(z)= za\ \mathcal{X}(z)+zb \ \mathcal{U}(z)\\
\mathcal{Y}(z)=c\ \mathcal{X}(z)+d \ \mathcal{U}(z).\\
\end{array}
\right.
\end{equation}
On $\mathbb{C}_I$ all
the objects commute, thus we have:
\begin{equation}\label{linqu4}
\left\{\begin{array}{l}
\mathcal{X}(z)=(1-za)^{-1} \ zb \  \mathcal{U}(z),\\
\mathcal{Y}(z)= c\ \mathcal{X}(z)+d \ \mathcal{U}(z).\\
\end{array}
\right.
\end{equation}

All the functions in the complex variable $z$ involved in system (\ref{linqu4}) are holomorphic on the plane $\mathbb{C}_I$,
thus we can use the
extension operator (\ref{ext}) to obtain a function in the quaternionic variable $p$ which is slice hyperholomorphic and since
$$
{\rm ext}(f(z)g(z))={\rm ext}(f(z)\star g(z))= f(p)\star g(p)
$$
we get
\begin{equation}\label{linqu41}
\left\{\begin{array}{l}
\mathcal{X}(p)= (1-pa)^{-\star}\star( pb)\star\mathcal{U}(p),\\
\mathcal{Y}(p)= c\ \star \ \mathcal{X}(p)+ d\ \star\ \mathcal{U}(p)\\
\end{array}
\right.
\end{equation}
from which
\begin{equation}\label{linqu42}
\left\{\begin{array}{l}
\mathcal{X}(p)= (1-pa)^{-\star}\star( pb)\star\mathcal{U}(p),\\
\mathcal{Y}(p)= c\ \star(1-pa)^{-\star}\star( pb)\star\mathcal{U}(p) + d\ \star\ \mathcal{U}(p).\\
\end{array}
\right.
\end{equation}
Since we have defined $\mathcal{H}(p)$ as
$$
\mathcal{H}(p):=  \mathcal{Y}(p)\star(\mathcal{U}(p))^{-\star}
$$
we obtain
$$
\mathcal{H}(p)=\Big(c\ \star(1-pa)^{-\star}\star( pb)\star\mathcal{U}(p) + d\
\star\ \mathcal{U}(p)\Big)\star (\mathcal{U}(p))^{-\star}
$$
$$
=c\ \star(1-pa)^{-\star}\star( pb) + d.
$$
This function is slice hyperholomorphic in $p$ with values in the quaternions so it makes sense to consider $a$, $b$, $c$ and $d$ as quaternions.
In order to get a matrix valued function  it is sufficient to
replace $a$, $b$, $c$ and $d$, respectively,  with the matrices $A$, $B$, $C$, $D$,  of appropriate size and with quaternionic entries.
We obtain  the quaternionic transfer function
$$
\mathcal{H}(p)=D+C\ \star(I-pA)^{-\star}\star( pB) =D+pC\ \star(I-pA)^{-\star} B
$$
which is a slice regular function.
Finally, using (\ref{eq:oberkampf_ligne_5}), we obtain the equality (\ref{reso}).
\end{proof}
\begin{Rk}{\rm
Using the explicit expression
for $(1-pA)^{-\star}$
we can also write
$$
\mathcal{H}(p)=D+C\ \star(I-\bar p A)(|p|^2A^2-2{\rm Re}(p) A
+I)^{-1}\star( pB),
$$
and if we  remove the $\star$-product we get $$\mathcal{H}(p)= D+(pC-|p|^2CA)(|p|^2A^2-2{\rm Re}(p) A
+I)^{-1}B.$$
}
\end{Rk}

\section{Rational functions}
\setcounter{equation}{0}
\label{sec:real}
Motivated by the discussion in the previous section, we now give
various equivalent definitions of rationality. We first consider
the case of $\mathbb H^{M\times N}$-valued functions of a real
variable, and give a definition and some preliminary lemmas.

\begin{Dn}
A function $f(x)$ of the real variable $x$ and with values in
$\mathbb H^{M\times N}$ will be said to be rational if it can be
written as
\begin{equation}
\label{eqrat}
f(x)=\frac{M(x)}{m(x)},
\end{equation}
where $M$ is a $\mathbb H^{M\times N}$-valued polynomial and
$m\in\mathbb R[x]$.
\end{Dn}
Clearly sums and products of rational functions of appropriate sizes stay rational. The case of the
inverse is considered in the next lemma.

\begin{La}
\label{la_real_rat}
Let $f$ be a rational function from $\mathbb R$ into $\mathbb H^{N\times N}$, and assume that
$f(x_0)$ is invertible for some $x_0\in\mathbb R$. Then, $f(x)$ is invertible for all $x\in\mathbb R$, at the possible
exception of a finite number of values, and $f^{-1}$ is a rational function.
\end{La}

\begin{proof} We proceed by induction. When $N=1$, a polynomial in $\mathbb{H}[x]$ can be seen as a matrix polynomial
in $\mathbb R^{4\times 4}[x]$, and this latter is invertible for all but at most a finite number of real values of $x$.
Let now $r(x)=\frac{t_1(x)}{t(x)}$, where $t_1\in\mathbb H[x]$ and $t\in\mathbb R[x]$. Assume $t\not\equiv 0$. Then,
since $x$ is real, we may write
\[
t^{-1}(x)=\frac{\overline{t(x)}}{t(x)\overline{t(x)}}=\frac{\overline{t}(x)}{(t\overline{t})(x)}.
\]
This ends the proof for $N=1$ since $(t\overline{t})(x)\in\mathbb R[x]$.\\

Assume now the induction proved at rank $N$ and let $R$ be a $\mathbb H^{(N+1)\times (N+1)}$-valued rational function,
invertible for at least one $x_0\in\mathbb R$. Then identifying $\mathbb H^{(N+1)\times (N+1)}$ with
$\mathbb R^{4(N+1)\times 4(N+1)}$ we see that $R$ is invertible for all, but at most a finite number
of real values of $x$
\[
R(x)=\begin{pmatrix}a(x)&b(x)\\ c(x)&d(x)\end{pmatrix}.
\]
Without loss of generality, we can assume that $a(x)$ is not
identically equal to $0$ (otherwise, multiply $R$ on the left or
on the right by a permutation matrix; this does not change the
property of $R$ or of $R^{-1}$ being rational). We write (see for
instance \cite[(0.3), p. 3]{Dym_CBMS})
\[
\begin{split}
\begin{pmatrix}a(x)&b(x)\\ c(x)&d(x)\end{pmatrix}=
\begin{pmatrix}1&0\\
c(x)a(x)^{-1}&I_n\end{pmatrix}\times\\
&\hspace{-1cm}\times
\begin{pmatrix}
a(x)&0\\0&d(x)-c(x)a(x)^{-1}b(x)\end{pmatrix}
\begin{pmatrix}1&a(x)^{-1}b(x)\\
0&I_n\end{pmatrix},
\end{split}
\]
and so $d(x)-c(x)a(x)^{-1}b(x)$ is invertible for all, but at most a finite number of, values $x\in\mathbb R$.
We have

\begin{equation}
\label{eq:R-1}
\begin{split}
R^{-1}(x)=
\begin{pmatrix}1&-a(x)^{-1}b(x)\\
0&I_n\end{pmatrix}\times\\
&\hspace{-5cm}\times\begin{pmatrix}
a(x)^{-1}&0\\0&(d(x)-c(x)a(x)^{-1}b(x))^{-1}
\end{pmatrix}
\begin{pmatrix}1&0\\
-c(x)a(x)^{-1}&I_n\end{pmatrix}.
\end{split}
\end{equation}
The induction hypothesis at rank $N$ implies that $(d(x)-c(x)a(x)^{-1}b(x))^{-1}$ is rational, and so is
$R^{-1}$, as seen from \eqref{eq:R-1}.
\end{proof}

\begin{La}
Consider a polynomial $M(p)\in \mathbb H^{N\times N}[p]$:
\[
M(p)=\sum_{j=0}^J p^jM_j.
\]
 Then
\[
M(p)=D+pC\star(I-pA)^{-\star}B,
\]
where $D=M_0$,
\[
A=\begin{pmatrix}0_{N}&I_{N}&0_{N}&\cdots\\
0_{N}&0_{N}&I_{N}&0_{N}&\cdots\\
&\vdots & & &\vdots\\
0_{N}&\cdots&\cdots&0_{N}&I_{N}\\
0_{N}&0_{N}&\cdots&0_{N}&0_{N}
\end{pmatrix},
\]
\[
B=\begin{pmatrix}0_{N}\\ 0_{N}\\ \vdots \\ I_{N}
\end{pmatrix},\quad
C=\begin{pmatrix}M_J&M_{J-1}&\cdots&M_1\end{pmatrix}.
\]
\end{La}
\begin{proof}
The equality easily follows from the fact that
$$
(I-pA)^{-\star}=\begin{pmatrix}I_{N
}&pI_{N}&p^2I_{N}&\cdots&p^JI_{N}\\
0_{N}&I_{N}&pI_{N}&\cdots&p^{J-1}I_{N}\\
&\vdots & & &\vdots\\
0_{N}&\cdots&\cdots&I_{N}&pI_{N}\\
0_{N}&0_{N}&\cdots&0_{N}&I_{N}
\end{pmatrix}.
$$
\end{proof}
Assume now that the $\mathbb H^{N\times N}$-valued function $f(p)$ can be written as
\begin{equation}
\label{eq:real}
 f(p)=D+p C\star(I-pA)^{-\star}B,
\end{equation}
where $A,B,C$ and $D$ are matrices with entries in $\mathbb H$
and of
appropriate sizes, and that $D$ is invertible. Then,
\begin{equation}
f(p)^{-1}=D^{-1}-pD^{-1}C\star(I-p(A-BD^{-1}C))^{-\star}B
\end{equation}
We also recall the following formula. See for instance \cite{bgk1}
for the case of rational functions of a complex variable. The
proof is as in the classical case, and will be omitted.

\begin{La}
Let
\[
 f_j(p)=D_j+pC_j\star (I_{N_j}-pA_j)^{-\star}B_j,\quad j=1,2
\]
be two functions admitting realizations of the form \eqref{eq:real}, and respectively $\mathbb H^{M\times N}$ and
$\mathbb H^{N\times R}$-valued. Then the $\mathbb H^{M\times R}$-valued function $ f_1\star  f_2$ can be written in
the form \eqref{eq:real}, with $D=D_1D_2$ and
\[
A=\begin{pmatrix}A_1&B_1C_2\\0&A_2\end{pmatrix},\quad B=\begin{pmatrix}B_1D_2\\ B_2\end{pmatrix},\quad
C=\begin{pmatrix}C_1&D_1C_2\end{pmatrix}.
\]
\label{realmult}
\end{La}
The proof of the following Lemma is immediate.

\begin{La}
Let
\[
 f_j(p)=D_j+pC_j\star(I_{N_j}-pA_j)^{-\star}B_j,\quad j=1,2
\]
be two functions admitting realizations of the form \eqref{eq:real}, and respectively $\mathbb H^{M\times N}$ and
$\mathbb H^{M\times R}$-valued. Then the $\mathbb H^{M\times (N+R)}$-valued function
$
\begin{pmatrix} f_1& f_2\end{pmatrix}
$
can be written in
the form \eqref{eq:real}, with $D=\begin{pmatrix}D_1&D_2\end{pmatrix}$ and
\[
A=\begin{pmatrix}A_1&0\\0&A_2\end{pmatrix},\quad B=\begin{pmatrix}B_1&0\\0& B_2\end{pmatrix},\quad
C=\begin{pmatrix}C_1&C_2\end{pmatrix}.
\]
\label{realmult1}
\end{La}

In the following theorem,
\[
 f(p)=\sum_{n=0}^\infty {\rm diag}\underbrace{(p^n,
p^n,\ldots , p^n)}_{M   \,\,{\rm times}}f_n,\]
which we will denote
for short by
\[
 f(p)=\sum_{n=0}^\infty p^nf_n,\]
is a series of the quaternionic variable $p$, with coefficients in
 $\mathbb H^{M\times N}$, converging in a neighborhood
of the origin.   We say that $ f$ is rational if any of the five conditions listed in the theorem holds.
\begin{Tm}
Let $ f$ be a $\mathbb H^{M\times N}$-valued function,
slice hyperholomorphic in a neighborhood of the origin. Then, the
following are equivalent:\\
$(1)$  There
exist matrices $A,B$ and $C$, of appropriate dimensions, and such
that
\begin{equation}
\label{goh}
f_n=CA^{n-1}B,\quad n=1,2,\ldots\qquad f_0=D.
\end{equation}
$(2)$ $ f$ can be written as \eqref{eq:real}
\[
 f(p)=D+p C\star(I-pA)^{-\star}B.
\]
$(3)$ The right linear span $\mathcal M( f)$ of the columns
of the functions $R_0f,R_0^2f,\ldots$
is a finite dimensional right quaternionic Hilbert space.\\
$(4)$ The function $ f(x)$ is a rational function from
$\mathbb R$ into $\mathbb H^{M\times N}$.\\
$(5)$ The entries of $ f$ are of the form $P\star
Q^{-\star}$, where $P$ and $Q$ are slice holomorphic polynomials,
and $Q(0)\not=0$.
\end{Tm}
\begin{proof} In the proof we adapt well known arguments from the
theory of matrix-valued rational functions to the present case.
For similar arguments in the classical case, see for instance \cite{bgk1}.\\

Assume that $(1)$ is in force, and set $D=s_0$. Then, for $p$
such that $|p|\cdot\|A\|<1$, the series $\sum_{n=0}^\infty p^n
A^n$ converges in ${\mathbb H}^{N\times N}$, and
\[
\sum_{n=0}^\infty p^n A^n=(I-pA)^{-\star},\ \ {\rm where} \ \ I=I_N,
\]
and we can write
\[
\sum_{n=0}^\infty p^nf_n=D+p C\star(I-pA)^{-\star} B,
\]
so that $(2)$ is in force.\\

Assume now $(2)$. Then,
\[
R_0^n f = C\star(I-pA)^{-1} A^{n-1}B,\quad n=1,2,\ldots,
\]
%\textcolor{red}{It was \[ R_0^n f = C(I-pA)^{-1}\star
%A^{n-1}B,\quad n=1,2,\ldots,
%\]}
so that $\mathcal M( f)$ is spanned by the columns of $D$
and the columns of the function $C\star(I-pA)^{-\star}$, and is
in particular finite dimensional, so that $(3)$ is in force.\\

Assume now $(3)$. Then there exists an integer $m_0\in\mathbb N$
such that for every $m\in\mathbb N$ and $v\in\mathbb H^q$, there
exist vectors $u_1,\ldots, u_{m_0}$ such that
\begin{equation}
\label{eq:generate}
R_0^{m_0} f v=\sum_{m=1}^{m_0}R_0^m
f u_m.
\end{equation}
Of course, the $u_j$ need not be unique.
Let $ E$ denote
the $\mathbb H^{p\times m_0 q}$-valued slice hyperholomorphic function
\[
 E=\begin{pmatrix}R_0 f&R_0^2
f&\ldots&R_0^{m_0} f\end{pmatrix}.
\]
Then, in view of \eqref{eq:generate}, there exists a matrix
$A\in\mathbb H^{m_0q\times m_0q}$ such that
\[
R_0 E= E A,
\]
so that
\[
 E(p)- E(0)=p\star  E(p) A= E(p)\star
p A,
\]
and so
\[
 E(p)= E(0)\star(I-p A)^{-\star}
\]
and
\[
(R_0  f)(p)= E(p)\begin{pmatrix}I_q\\ 0\\
\vdots\\0\end{pmatrix}.
\]
Thus we have
\[
f(p)-f(0)=p E(0)\star(I-p A)^{-\star}\begin{pmatrix}I_q\\ 0\\
\vdots\\0\end{pmatrix}.
\]
It follows that $ f$ is of the form \eqref{eq:real}, and
so $(2)$ holds. Since $(2)$ implies trivially $(1)$, we have
shown that the first three claims are equivalent. We now turn to
the equivalence with the other two characterizations.\\

Assume that $(3)$ holds. Then, the restriction of \eqref{eq:real} to $p=x\in\mathbb R$ gives
\[
 f(x)=D+xC(I-xA)^{-1}B.
\]
and thus $ f$ is a rational function of $x$ in view of Lemma \ref{la_real_rat}, that is $(4)$ is in force. Assume now
$(4)$, and let $ f(x)=\frac{M(x)}{m(x)}$ as in \eqref{eqrat}. Since $ f$ is assumed defined the origin, we may
assume that $m(0)\not=0$, and $ f(x)$ is the restriction to the real line of the slice hyperholomorphic
function
\[
 f(p)=M(p)\star m(p)^{-\star} =m(p)^{-\star} \star M(p).
\]
Note that since $m$ has real coefficients we also have $m(p)^{-\star} \star M(p)=m(p)^{-1}  M(p)$.
So $(5)$ holds. Assume now $(5)$ in force. Then, in view of Lemma \ref{realmult}, each of the entries of $ f$ admits
a realization, and the result follows since the property of realization stays under concatenation. Use Lemma
\ref{realmult1} and its counterpart for columns rather than rows.
\end{proof}

We note that a short introduction to realization theory in forms of a sequence of
exercises can be found in \cite[pp. 328-329]{alpay_book}.

\section{Positive definite functions and reproducing kernels}
\setcounter{equation}{0}
\label{sec:positive}
The $\mathbb H$-valued function $k(u,v)$ defined for $u,v$ in some set $\Omega$ (in the present paper, $\Omega$ will be most of the time the unit ball of $\mathbb H$) is called positive
definite if it is Hermitian:
\begin{equation}
\label{herm}
k(u,v)=\overline{k(v,u)},\quad\forall u,v\in\Omega,
\end{equation}

and if
for every $N\in\mathbb N$, every $u_1,\ldots , u_N\in \Omega$ and $c_1,\ldots, c_N\in\mathbb H$ it holds that
\[
\sum_{\ell,j=1}^N\overline{c_\ell}k(u_\ell, u_j)c_j\ge 0.
\]
(Note that the above sum is a real number in view of \eqref{herm}). As in the complex-valued case, associated to $k$ is a uniquely defined reproducing kernel quaternionic (right)-Hilbert space $\mathscr H(k)$, with reproducing kernel $k(u,v)$, meaning that:\\
$(1)$ The function $u\mapsto k(u,v)c$ belongs to $\mathscr H(k)$ for every choice of $v\in \Omega$ and $c\in\mathbb H$, and\\
$(2)$ it holds that
\[
\overline{c}f(v)=\langle f(\cdot), k(\cdot, v)c\rangle_{\mathscr H(k)}
\]
for every choice of $f\in\mathscr H(k)$ and of $u\in\Omega$ and $c\in\mathbb H$.
See for instance \cite[Theorem 8.4, p. 456]{as3} for a proof and
\cite[Proposition 9.4, p. 458]{as3} for a proof of a result that, with similar techniques, proves also the following lemma.

\begin{La}
\label{la:element}
Let $\mathscr H(k)$ be a right reproducing kernel quaternionic Hilbert space of $\mathbb H$-valued functions defined on a set $\Omega$, with reproducing kernel $k(u,v)$. Then, the
function $f$ belongs to $\mathscr H(k)$ if and only if there exists $M>0$ such that the kernel
\[
k(u,v)-\frac{f(u)\overline{f(v)}}{M^2}\ge 0.
\]
The smallest such $M$ is equal to the norm of $f$ in $\mathscr H(k)$.
\end{La}

Reproducing kernels are a main tool in
operator theory. Reproducing kernels of the form $c(z\bar w)$,
$z,w\in\mathbb{C}$, where $c(t)=\sum_{n=0}^{+\infty} c_n t^n$,
$c_n\geq 0$, for all $n\in\mathbb{N}$ is an analytic function is
a neighborhood of the origin are an important case. In the case
of Hardy space, Bergman space, and Dirichlet space the functions
$c(t)$ are given by
$$
c(t)=\frac{1}{1-t}\qquad c(t)=\frac{1}{(1-t)^2} \qquad
c(t)=-\ln(1-t).
$$
We will consider the generalization of these kernels to the case
in which the variables considered are quaternions, i.e. we will
consider
\begin{equation}\label{series}
c(p\bar q)=\sum_{n=0}^{+\infty} c_n p^n\bar q^n,\qquad
q,p\in\mathbb{H}.
\end{equation}
Note that due to the noncommutative nature of quaternions, the
series in (\ref{series}) does not coincide with
$\sum_{n=0}^{+\infty} c_n (p\bar q)^n$.

\begin{Tm}
Let $A\subset\mathbb N_0$ and let $(c_n)_{n\in A}$ be a sequence of strictly
positive numbers, and assume that
\[
R=\frac{1}{\limsup_{n\in A} c_n^{1/n}}>0.
\]
 Then, the function
\[
k(p,q)=\sum_{n\in A}c_np^n\overline{q}^n
\]
is positive definite in the ball $|p|<R$. The associated reproducing kernel
Hilbert space consists of the functions
\[
f(p)=\sum_{n\in A}p^nf_n,
\]
with coefficients $f_n\in\mathbb H$ such that
\[
\sum_{n\in A}\frac{|f_n|^2}{c_n}<\infty.
\]
\end{Tm}

%The proof is already known in the case of the Bergman kernel, see ....., in which we have also proved that $K(q,p)=....$.

In the case of the Hardy space we have
\begin{Pn}
The sum of the series $\sum_{n=0}^{+\infty} p^n\bar q^n$ is the function $k(p,q)$ given by
\begin{equation}\label{kernel}
k(p,q)=(1-2{\rm Re}(q) p+|q|^2p^2)^{-1}(1-pq)=(1-\bar p \bar q) (1-2{\rm Re}(p) \bar q+|p|^2\bar q^2)^{-1}.
\end{equation}
The kernel $k(p,q)$ is defined for $p\not\in[q^{-1}]$ or, equivalently, for $q\not\in[p^{-1}]$.
Moreover:
\begin{itemize}
\item[a)] $k(p,q)$ is left slice regular in $p$ and right slice regular in $q$;
\item[b)] $\overline{k(p,q)}=k(q,p)$.
\end{itemize}
\end{Pn}
 \begin{proof}
 Consider the function
 $$
 k_q(z)=\frac{1}{1-z\bar q}, \qquad z=x+Iy.
 $$
 The proof  of the first equality is an application of the representation formula (\ref{strutturaquat}).
 Consider now the function
 $$
 k_p(w)=\frac{1}{1-p\bar w};
 $$
 the representation formula (\ref{distribution}) gives
 $k(p,q)=(1-\bar p \bar q) (1-2{\rm Re}(p) \bar q+|p|^2\bar q^2)^{-1}$. The function
$(1-\bar p \bar q) (1-2{\rm Re}(p) \bar q+|p|^2\bar q^2)^{-1}$ is right slice regular in the variable $q$ by construction and it is slice regular in the variable $p$ in its domain of definition, since it is the product of a slice regular function and a polynomial with real coefficients. By the identity principle it coincides with the first expression which is slice regular in $p$ by construction. Assertion a) follows.  Point b)  follows from the chain of equalities
$$\overline{k(q,p)}=\overline{[(1-\bar p \bar q) (1-2{\rm Re}(p) \bar q+|p|^2\bar q^2)^{-1}]}
$$
$$
=
(1-2{\rm Re}(p)  q+|p|^2\bar q^2)^{-1}(1- q p)=k(q,p).$$
\end{proof}

We note (see for instance \cite{donoghue,MR2002b:47144}):

\begin{Tm}
 The
following are equivalent:
\\
$(1)$ The function $s$ is analytic and contractive in the open
unit disk.\\
$(2)$  The function $s$ is {\sl defined} in $\mathbb D$ and the
operator of multiplication by $s$ is a contraction from the Hardy
space $ H_2(\mathbb D)$ into itself.\\
$(3)$ The function $s$ is {\sl defined} in $\mathbb D$ and the
kernel
\[
k_s(z,w)=\frac{1-s(z)\overline{s(w)}}{1-z\overline{w}}
\]
is positive in the open unit disk.
\label{Tm:donoghue}
\end{Tm}

We prove in Theorem \ref{donoghueh} the counterpart of the above theorem in the
quaternionic setting.\\

 Note that $k_s$ can be
rewritten as
\[
k_s(z,w)=\sum_{n=0}^\infty
z^n(1-s(z)\overline{s(w)})\overline{w}^n.
\]
This is the form which will be used in the quaternionic setting.
\
See \eqref{ksh} in Section \ref{sec:schur_mul}.\\

The reproducing kernel Hilbert space with reproducing kernel
$k_s(z,w)$ was first introduced and studied by de Branges and
Rovnyak; see \cite{dbr2, dbr1}. This space will be denoted by
$\mathscr H(s)$. It is equal to the operator range ${\rm
ran}~\sqrt{I-M_sM_s^*}$, where $M_s$ denotes the operator of
multiplication by $s$ from $ H_2(\mathbb D)$ into itself,
that is, $M_s$ is the commutative version of the operator defined
by \eqref{ms_h} in the
following section.\\

The space $\mathscr H(s)$ is the state space for a coisometric
realization of $s$. Furthermore, $s$ is a Schur function of and
only if the function $s^\sharp$ defined by
\begin{equation}
\label{sharp}
s^\sharp(z)=\overline{s(\overline{z})}
\end{equation}
is a Schur function. The space $\mathscr H(s^\sharp)$ is the state
space for an isometric realization of $s$. A unitary realization
for $s$ is given in terms of the reproducing kernel Hilbert space
with reproducing kernel
\begin{equation}
\label{eq:DS}
D_s(z,w)=\begin{pmatrix}
k_s(z,w)&\displaystyle\frac{s(z)-s(\overline{w})}{z-\overline{w}}\\
\displaystyle\frac{s^\sharp(z)-s^\sharp(\overline{w})}{z-\overline{w}}&k_{s^\sharp}(z,w)
\end{pmatrix}.
\end{equation}

For more information on the Schur algorithm, see for instance the
papers \cite{kailath} and the books
\cite{MR2002b:47144,MR94b:47022,C-schur}.

\section{Schur multipliers and Schwarz' lemma}
\setcounter{equation}{0}
\label{sec:schur_mul}
%\section{Linear systems in the quaternionic case}
%\setcounter{equation}{0}
We will define a finite energy quaternionic signal as a sequence $x=(x_0,x_1,\ldots)$ of quaternions such that
\[
\|{x}\|^2\stackrel{\rm def.}{=}\sum_{n=0}^\infty |x_n|^2<\infty,
\]
that is, an element of $\ell_2(\mathbb N_0,\mathbb H)$.  Note that
$\ell_2(\mathbb N_0,\mathbb H)$ is a quaternionic right Hilbert
space when endowed with the $\mathbb H$-valued inner product
\[
[x,  y]=\sum_{n=0}^\infty \overline{y_n}x_n,\quad{\rm where}\quad
 y=(y_0,y_1,\ldots).
\]
To $x\in\ell_2(\mathbb N_0,\mathbb H)$ we associate the function
\[
f_{x}(p)=\sum_{n=0}^\infty p^nx_n
\]
of the quaternionic variable $p$. We denote by $\mathbf{H}_2$ the right quaternionic vector space of
such power series,
endowed with the inner product
\[
[f_x,f_y]=\sum_{n=0}^\infty\overline{ y_n}x_n.
\]
To ease the notation, we will use the
notation $x(p)$ rather than $f_{x}(p)$:
\begin{equation}
x(p)=\sum_{n=0}^\infty p^nx_n.
\label{eq:iso}
\end{equation}

\begin{Dn}
We set $\ell(\mathbb N_0, \mathbb H)$ denote the set of sequences
of quaternions, indexed by $\mathbb N_0$, and  $\ell_0(\mathbb
N_0, \mathbb H)$ denote the set of sequences of $\ell(\mathbb N_0,
\mathbb H)$  for which at most a finite number of elements are
non-zeros.
\label{def:ells}
\end{Dn}

The slice regular product of $
x\in\ell_2(\mathbb N_0, \mathbb H)$ and $
y\in\ell_0(\mathbb N_0, \mathbb H)$ is defined by
\[
(x\star y)_n=\sum_{m=0}^n x_my_{n-m},\quad n=0,1,\ldots
\]
Note that
\[
x\star 1=x,\quad{\rm where}\quad  1=(1,0,0,\ldots).
\]

We now state the counterpart of Theorem \ref{Tm:donoghue} in the quaternionic
setting:

\begin{Tm}
Let $ s\in\ell(\mathbb N_0, \mathbb H)$ (see Definition \ref{def:ells}). Then, the
following are equivalent:\\
$(1)$ The map $x\mapsto  s\star x$, first
defined for $x\in \ell_0(\mathbb N_0, \mathbb H)$,
extends to a contraction from $\ell_2(\mathbb N_0, \mathbb H)$
into itself.\\
$(2)$ The series $ s(p)=\sum_{n=0}^\infty p^ns_n$
converges in $\mathbb B$ and the operator
\begin{equation}
\label{ms_h} (M_{ s}x)(p)=\sum_{n=0}^\infty
p^n s(p)x_n
\end{equation}
is a contraction from the quaternionic Hardy space $\mathbf H_2(\mathbb B)$ into itself, where $
x(p)=\sum_{n=0}^\infty p^nx_n\in \mathbf H_2(\mathbb B)$.\\
$(3)$ The series $ s(p)=\sum_{n=0}^\infty p^ns_n$
converges in $\mathbb B$, and the function
\begin{equation}
\label{ksh} k_{ s}(p,q)=\sum_{n=0}^\infty p^n(1-{
s}(p)\overline{{ s}(q)})\overline{q}^n\\
= (1-s(p)\overline{s(q)}) \star (1-p\overline{q})^{-\star}
\end{equation}
is positive definite in the unit ball $\mathbb B$ of $\mathbb H$.
\label{donoghueh}
\end{Tm}

\begin{proof} We first assume that $(1)$ is in force.
Since $ s\star 1= s$ we note that in fact
$ s\in\ell_2(\mathbb N_0, \mathbb H)$. The power series
$ s(p)$ makes thus sense for every $p$ in the open unit
ball of $\mathbb H$. We denote by $M_{ s}$ the image of the
operator
 $x\mapsto  s\star x$ under the map
 \eqref{eq:iso}. Then, for $ x\in \ell_0(\mathbb N_0, \mathbb H)$
 we have
\begin{equation}
\label{eq:Ms}
 \sum_{n=0}^\infty p^n s(p)x_n=\sum_{n=0}^\infty p^n
(\sum_{m=0}^\infty p^m s_m) x_n=\sum_{n=0}^\infty p^n (
s\star  x)_n,
\end{equation}
where the sum is finite since $ x\in \ell_0(\mathbb N_0,
\mathbb H)$. Since in a reproducing kernel Hilbert space,
convergence in norm implies pointwise convergence, we see that
the operator $M_{s}$ is given by \eqref{eq:Ms} for every $ x\in
\ell_2(\mathbb N_0, \mathbb H)$, and hence $(2)$ is in force.\\

Assume now that $(2)$ holds. Consider the kernel $k(p,q)$  in (\ref{kernel}), which on $|pq|<1$
 admits the power series expansion
$$
k(p,q)=\sum_{n=0}^\infty p^n\overline{q}^n.
$$
For sake of simplicity we will write $k_q(p)$ instead of
$k(p,q)$. Let $q_1,q_2\in\mathbb B$ and $c_1,c_2\in\mathbb H$
 and compute
\[
\begin{split}
\overline{c_2}\left(M_{ s}^*k_{q_1}c_1\right)(q_2)&=
[
M_{ s}^*k_{q_1}c_1\, ,\, k_{q_2}c_2]_{\mathbf H_2(\mathbb B)}\\
&=[
k_{q_1}c_1
\, ,\, M_{ s}(k_{q_2}c_2)]_{\mathbf H_2(\mathbb B)}\\
&=
[
k_{q_1}c_1\, ,\, \sum_{n=0}^\infty p^n s(p)\overline{q_2}^nc_2]_{\mathbf H_2(\mathbb B)}\\
&=
\overline{  [\sum_{n=0}^\infty p^n s(p)\overline{q_2}^nc_2 \, ,\,   k_{q_1}c_1 ]_{\mathbf H_2(\mathbb B)}}\\
&=\overline{c_2}\left(\sum_{n=0}^\infty q_2^n\overline{p_1}\overline{q_1}^nc_1\right).
\end{split}
\]
Thus, the adjoint of $M_{ s}$ is given by the formula:
\[
(M_{ s}^*(k(q))(p)=\sum_{n=0}^\infty
p^n\overline{s(q)^n}\overline{q}^n,
\]
and $(3)$ is just a rewriting that the operator $I-M_{
s}M_{ s}^*$ is non negative.\\

Assume now that $(3)$ holds. We consider the relation $R_{
s}$ (that is, the linear subspace) in $\mathbf H_2(\mathbb B)\times \mathbf
H_2(\mathbb B)$ spanned by the pairs

\[
\left(\sum_{n=0}^\infty p^n\overline{q}^nc\, ,\, \sum_{n=0}^\infty
p^n{ s}(p)\overline{{ s}(q)}\overline{q}^nc\right),
\]
where $p$ runs in $\mathbb B$ and $c$ in $\mathbb H$.
The domain of $R_{ s}$ is dense, and the positivity of the
kernel implies that $R_{ s}$ is a contraction, meaning that if
$(f,g)\in R_{ s}$, then
\[
\|g\|_{\mathbf H_2(\mathbb B)}\le \|f\|_{\mathbf H_2(\mathbb B)}.
\]
It follows that $R_{ s}$ extends to the graph of an everywhere defined
contraction, which we will denote by $T$. We now compute $T^*$:
\[
\begin{split}
\overline{c_2}(T^*k_{q_1}c_1)(q_1)&=
[T^*k_{q_1}c_1\, ,\, k_{q_2} c_2]_{\mathbf H_2(\mathbb B)}\\
&=[k_{q_1}c_1\, ,\, \sum_{n=0}^\infty
p^n\overline{{ s}(q_2)}\overline{q_2}^nc_2]_{\mathbf H_2(\mathbb B)}\\
&=\overline{c_2}\left(\sum_{n=0}^\infty q_2
s(q_2)\overline{q_1}^nc_1\right)\\
&=\overline{c_2}\left(M_{ s}(k_{q_2}c_1)\right)(q_1).
\end{split}
\]
Thus $T^*=M_{ s}$. We obtain $(1)$ by looking at the
operator $M_{ s}$ on each $p^n$.
\end{proof}

We note that the kernel $k_{ s}$ can also be written as

\[
\begin{split}
k_{ s} (p,q)&=(1-2{\rm Re}(q) p+|q|^2p^2)^{-1}(1-pq)\star
(1-s(p)\overline{s(q)})
\\
&=(1-2{\rm Re}(q)
p+|q|^2p^2)^{-1}(1-pq-s(p)\overline{s(q)}+ps(p)q\overline{s(q)}).
\end{split}
\]

The operator $M_{ s}$ is the non-commutative version of the
operator of multiplication by ${s}$. We note the following
properties of the operators $M_{ s}$:

\begin{Pn}
Let $ s_1$, $ s_2$ and $ s$ be Schur multipliers depending on  the quaternionic variable $p$.
Then:
\begin{eqnarray}
M_{{s}_1}M_{{s}_2}&=&M_{{s}_1\star {s}_2},\\
M_{{s}}M_{p}&=&M_{p}M_{s}=M_{p{s}}.
        \end{eqnarray}
\end{Pn}

\begin{proof}
The first equality follows from the associativity of
the $\star$-product. The second property is a consequence of Remark \ref{gcomf} while the last one follows from the previous ones.
\end{proof}

The proof of the following proposition follows the case of analytic functions
in the open unit disk (see \cite[Lemma 1, p. 301]{MR0140931}
for instance for the latter), and is
given for completeness. In the case of analytic functions
in a half-plane, the corresponding result is more difficult to prove,
and is called the Bochner-Chandrasekharan theorem. See
\cite{boch_chan, weiss_MCSS}.

\begin{Rk}{\rm We now discuss a feature which arises because we are working in
a non-commutative setting.
When considering a vector space on $\mathbb{H}$,  the set of
right linear operators acting on $V$ has a structure of vector
space over $\mathbb{H}$ only if $V$ is a two sided vector space,
in fact right linear operators form a left vector space, see
Section 2. However, if in $V$ it is defined a multiplication
$\star$ among vectors compatible with the right vector space
structure, or at least, if  for some $w\in V$ the multiplication
$w\star v$ is an operation in $V$ for all $v\in V$, then this
multiplication allows to define right linear operators which form
a right linear vector space. Thus, in this case, it is not
necessary to require that $V$ is two sided vector space. Indeed,
we can define a right linear operator $M_w$:
$$
M_w(v):=w\star v.
$$
The linearity is obvious and the structure of right linear space on the operators is given by
$$
(M_w\lambda)(v):=w\lambda\star v, \qquad \lambda\in\mathbb{H}.
$$
In particular, if $V=\mathbf H_2(\mathbb
B)$ and $M$ is the multiplication on the left by $s(p)\in V$ (assuming that the multiplication is defined for every element in $V$) then we have $M_s(f)(p)=s(p)\star f(p)$ and if $\pi(p)=a_0+pa_1+\cdots +p^na_n$ then $M_s(\pi)(p)=s(p)\star \pi(p)=s(p)a_0+ps(p)a_1+\ldots+p^ns(p)a_n$.
}
\end{Rk}

\begin{Pn}
Let $T$ be a right-linear contraction from $\mathbf H_2(\mathbb
B)$ into itself, which commutes with $M_p$. Then $T=M_{
s}$ for some Schur multiplier ${ s}$.
\end{Pn}

\begin{proof}
Denote ${ s}(p)=(T1)(p)$. Then, for $a\in\mathbb
H\subset\mathbf H_2(\mathbb B)$ (in other words, in \eqref{eq:ta}
and \eqref{eq:ta2} below we view $a\in\mathbb H$ as a function in
$\mathbf H_2(\mathbb B)$ and, by the assumed right linearity, we
have $Ta=(T1)a$)
\begin{equation}
\label{eq:ta}
\left((TM_p)(a)\right)(p)=(M_p(Ta))(p)=((M_p(T1))a)(p)=p{
s}(p)a.
\end{equation}
More generally, an induction shows that
\begin{equation}
\label{eq:ta2}
\left((TM_{p^n})(a)\right)(p)=(M_{p^n}(Ta))(p)=p^n{ s}
(p)a,\quad n=2,3,\ldots
\end{equation}
Thus, for every polynomial
\[
\pi(p)=a_0+pa_1+\cdots +p^na_n,
\]
we have
\begin{equation}
\label{eq:Felix_Faure_Ligne_8}
(T\pi)(p)=\sum_{j=0}^n p^j{ s}(p)a_j.
\end{equation}
Let now $f\in\mathbf H_2(\mathbb B)$, with power series expansion
\[
f(p)=\sum_{j=0}^\infty p^ja_j,
\]
and let $f_n$ be the polynomial
\[
f_n(p)=\sum_{j=0}^n p^ja_j.
\]
We have $\lim_{n\rightarrow\infty}\|f-f_n\|_{{\mathbf H_2(\mathbb B)}}=0$, and so, by
continuity of $T$,
\[
\lim_{n\rightarrow\infty}\|Tf-Tf_n\|_{{\mathbf H_2(\mathbb B)}}=0
\]
In a reproducing kernel Hilbert space convergence in norm implies
pointwise convergence. Thus we have that, for every $p\in\mathbb B$,
\[
(Tf)(p)=\sum_{j=0}^\infty p^j{ s}(p)a_j,
\]
where we have used \eqref{eq:Felix_Faure_Ligne_8}. Thus, $T=M_{ s}$.
Finally, $ s$ is a Schur multiplier since $T$ is assumed contractive.
\end{proof}

The following result is a counterpart of Schwarz's lemma for Schur multipliers.

\begin{Tm}
\label{schwarz_slice}
Let ${ s}$ be a Schur multiplier, and assume that
$s(0)=0$. Set $ s(p)=p{s}^{(1)}(p)$.
Then ${s}^{(1)}$ is also a Schur multiplier.
\end{Tm}

\begin{proof} Let $\mathscr H(s)$ be the reproducing kernel
quaternionic Hilbert space with
reproducing kernel $k_{s}(p,q)$. Since $s(0)=0$ we have that
$1=k_{s}(p,0)\in\mathscr H({ s})$,  and
$\|1\|_{\mathscr H({ s})}^2=k_{ s}(0,0)=1$. Hence the kernel
\[
k_{ s}(p,q)-1
\]
is positive definite in $\mathbb B$; see Lemma \ref{la:element}. Thus
\[
\sum_{n=0}^\infty p^n\left(1-p{ s}^{(1)}(p)\overline{{ s}^{(1)}(q)}
\overline{q}\right)\overline{q}^n\ge 1,
\]
which can be rewritten as
\[
p\left(\sum_{n=0}^\infty  p^n\left(1-{ s}^{(1)}(p)\overline{
{ s}^{(1)}(q)}\right)\overline{q}^n\right)\overline{q}\ge 0.
\]
Thus,
\[
\sum_{n=0}^\infty  p^n(1-{ s}^{(1)}(p)\overline{{ s}^{(1)}(q)})
\overline{q}^n
\ge 0,
\]
and ${ s}^{(1)}$ is a Schur multiplier.
\end{proof}

We conclude this section with the following characterization of
the space $\mathscr H(s)$:

\begin{Tm}
Let $s$ be a Schur multiplier. Then the associated reproducing
kernel Hilbert space $\mathscr H(s)$ is equal to the operator
range ${\rm ran}~(I-M_sM_s^*)^{1/2}$ endowed with the norm
\[
\|(I-M_sM_s^*)^{1/2}u\|_{\mathscr H(s)}=\|(I-\pi)u\|_{{\mathbf
H_2(\mathbb B)}},
\]
where $\pi$ is the orthogonal projection on the kernel of
$(I-M_sM_s^*)^{1/2}$.
\end{Tm}

The proof is as in the classical case and will be omitted. See
\cite[Theorem 3.2 p. 16]{ad3} for instance and \cite[Theorem 8.3,
p. 119]{MR93b:47027} for a similar proof in the setting of upper
triangular operators. We refer to  \cite{fw} for more on operator
ranges and for their relations with de Branges Rovnyak spaces.

\section{Realizations of Schur multipliers}\label{sec:co}
\setcounter{equation}{0}
The purpose of this section is to prove the following result,
which is the counterpart in the quaternionic setting of the
realization result for Schur functions. We first recall the
following definition (see for instance \cite[p. 14]{adrs} in the
complex case).

\begin{Dn}
Let $\mathscr H$ be a right quaternionic Hilbert space. The pair
$(C,A)\in\mathbf L(\mathscr H,\mathbb H^p)\times \mathbf L
(\mathscr H)$ is called closely outer connected if
\begin{equation}
\bigvee_{n=0}^\infty {\rm ran}~A^{*n}C^*=\mathscr H.
\end{equation}
\end{Dn}

The following is the counterpart of \cite[Theorem 2.2.1, p.
49]{adrs} in the setting of slice holomorphic functions.

\begin{Tm}
\label{Vaugirard} Let $s$ be a function from the open unit ball
of $\mathbb H$ into $\mathbb H$. Then, $s$ is a Schur multiplier
if and only if there exist a right quaternionic Hilbert space
$\mathscr H$ and a coisometric operator
\begin{equation}
\label{eq:ABCD}
\begin{pmatrix}A&B\\ C&D\end{pmatrix}\,\,:\,\,\mathscr
H\oplus\mathbb H\longrightarrow \mathscr H\oplus\mathbb H
\end{equation}
such that ${s}$ can be written as a power series
\[
s(p)=\sum_{n=0}^\infty p^ns_n
\]
where
\begin{equation}
\label{eq:sn}
s_n=\begin{cases}D,\quad\hspace{1.2cm} n=0,\\
CA^{n-1}B1,\quad n=1,2,\ldots \
\end{cases}
\end{equation}
Assume the pair $(C,A)$ closely outer connected. Then, it is unique up
to an isometry of right quaternionic Hilbert spaces.
\end{Tm}

We note that the realization \eqref{eq:sn} is called closely
outer connected when the pair $(C,A)$ is closely outer
connected.\\

Using  formula \eqref{eq:oberkampf_ligne_5} we have for $\|pA\|<1$
\[
\sum_{n=1}^\infty  p^nCA^{n-1}B=C\star (\sum_{n=1}^\infty
p^nA^{n-1})B=C\star S^{-1}_R(p^{-1},A)B,
\]
and so
\begin{equation}
s(p)=D+ C\star S^{-1}_R(p^{-1},A)B=D+pC\star (I-pA)^{-*}B.
\end{equation}

The strategy of the proof is as follows.
Let $\s$ be a Schur multiplier, and let $\mathscr H(\s)$ be the
associated reproducing kernel quaternionic Hilbert space. As in
the classical case, we want to show that $\mathscr H(\s)$ is the
state space, in an appropriate sense in the present setting, of a
coisometric realization of $\s$. We use the same method as in
\cite{adrs} (see in particular p. 50 there),
suitably adapted to the non commutativity of
$\mathbb H$. One considers the reproducing kernel quaternionic
right Hilbert space with reproducing kernel $k_{ s}(p,q)$,
and define the relation $R$ (that is the right vector subspace of
$(\mathscr H(\s)\oplus\mathbb H)\times (\mathscr H(\s)\oplus\mathbb
H)$ defined as the right linear span of elements of the form
\[
\left\{ \begin{pmatrix} k_\s(p,q)\overline{q}u\\
\overline{q}v\end{pmatrix},
\begin{pmatrix}(k_\s(p,q)-k_\s(p,0))u+k_\s(p,0)\overline{q}v\\
(\overline{\s(q)}-\overline{\s(0)})u+
\overline{\s(0)}\overline{q}v\end{pmatrix}\right\}.
\]
We claim that the relation $R$ is densely defined and isometric.
It will follow that $R$ can be extended in a unique way to the
graph of an isometric operator from $\mathscr H(\s)\oplus\mathbb H$
into itself. This operator (or more precisely its adjoint) will
give the realization. We first prove some preliminary lemmas and then
give the proof of the theorem along the above lines.

\begin{La}
The relation $R$ is isometric and densely defined.
\end{La}
\begin{proof}
We want to prove that
\begin{equation}
\label{bastille}
\begin{split}
[ k_s(p,q_1)\overline{q_1}u_1,
k_s(p,q_2)\overline{q_2}u_2]_{\mathscr H(s)}+
[\overline{q_1}v_1,\overline{q_2}v_2]_{\mathbb H}&=\\
&\hspace{-8cm}= [
(k_s(p,q_1)-k_s(p,0))u_1+k_s(p,0)\overline{q_1}v_1,(k_s(p,q_2)-k_s(p,0))u_2+
k_s(p,0)\overline{q_2}v_2]_{\mathscr H(s)}+\\
&\hspace{-7.5cm}+
[(\overline{s(q_1)}-\overline{s(0)})u_1+\overline{s(0)}
\overline{q_1}v_1,(\overline{s(q_2)}-\overline{s(0)})u_2
+\overline{s(0)}\overline{q_2}v_2]_{\mathbb H}
\end{split}
\end{equation}
for all choices of $u_1,u_2,v_1,v_2\in\mathscr H$ and $q_1,q_2$ in
the open unit ball of $\mathbb H$. We rewrite this equality as
\[
\overline{u_2}M_1 u_1+\overline v_2M_2 v_1+\overline{v_2}M_3u_1
+\overline{u_2}M_4v_1=\overline{u_2}N_1 u_1+\overline v_2N_2 v_1+
\overline{v_2}N_3u_1
+\overline{u_2}N_4v_1,
\]
and we will show that $M_i=N_i$ for $i=1,2,3,4$. We write out
in details the case $i=1$, and only outline the other cases.\\

{\bf Checking $M_1=N_1$:}
Using the reproducing kernel property, we see
that
\[
M_1=q_2k_{ s}(q_2,q_1)\overline{q_1}.
\]
Still with this same property, we have
\[
\begin{split}
N_1&=k_{ s}(q_2,q_1)-k_{ s}(q_2,0)-
k_{ s}(0,q_1)+k_{ s}(0,0)+\\
&\hspace{5mm}+(s(q_2)-s(0))(\overline{s(q_1)}-\overline{s(0)})\\
&=k_{ s}(q_2,q_1)-(1-s(q_2)\overline{s(0)})-(1-
s(0)\overline{s(q_1)})+1-|s(0)|^2+\\
&\hspace{5mm}+
s(q_2)\overline{s(q_1)}-s(q_2)\overline{s(0)}-
s(0))\overline{s(q_1)}+|s(0)|^2
\end{split}
\]
and thus the to prove the equality $M_1=N_1$ is equivalent to check that
\[
k_s(q_2,q_1)-q_2k_s(q_2,q_1)\overline{q_1}=1-
s(q_2)\overline{s(q_1)},
\]
but this is a direct consequence of the definition of the kernel $k_s$.\\

{\bf Checking $M_2=N_2$:}  We now have $M_2=q_2\overline{q_1}$ and
\[
N_2=q_2k_s(0,0)\overline{q_1}+q_2|s(0)|^2\overline{q_1}=
q_2\overline{q_1}=M_2.
\]

We have $M_3=M_4=0$ and so we now need to check that $N_3=N_4=0$.\\

{\bf Checking $N_3=0$:} This amounts to verify that
\[
q_2(k_s(0,q_1)-k_s(0,0))+q_2(s(0)(\overline{s(q_1)}-\overline{s(0)})=0,
\]
but this is clear from the definition of the kernel since
\[
k_s(0,q_1)=1-s(0)\overline{s(q_1)}\quad{\rm and}\quad k_s(0,0)=1-|s(0)|^2.
\]

{\bf Checking $N_4=0$:} This amounts to verify that
\[
(k_s(0,0)-k_s(0,q_2))\overline{q_1}+(s(0)-s(q_2))\overline{s(0)}
\overline{q_1}=0,
\]
but this is also plain from the definition of the kernel $k_s$.
\end{proof}

\begin{La}
$R$ is the graph of a densely defined isometry. Let us denote by
\[
\begin{pmatrix}A&B\\ C&D\end{pmatrix}^*\,\,:\,\,\mathscr
H(s)\oplus\mathbb H\longrightarrow \mathscr H(s)\oplus\mathbb H
\]
its extension to all of $\mathscr H(s)\oplus\mathbb H$. Then
\begin{equation}
\label{eq:dbh}
\begin{split}
(Af)(p)&=\begin{cases}p^{-1}(f(p)-f(0)),\quad p\not =0\\
f_1,\quad\hspace{2.6cm}
 p=0,
\end{cases}\\
 (Bv)(p)&=\begin{cases}p^{-1}(s(p)-s(0))v,
 \quad p\not =0\\
s_1,\quad\hspace{2.6cm}p=0,
 \end{cases}\\
 Cf&=f(0),\\
 Dv&=s(0)v.
 \end{split}
\end{equation}
\end{La}

\begin{proof} We first check that $R$ is the graph of a densely defined
isometry. By definition, the domain of $R$ is the set of $F\in\mathscr H(s)
\oplus\mathbb H$ such that there exists $G\in
\mathscr H(s)
\oplus\mathbb H$ such that $(F,G)\in R$. It is therefore dense by construction
since the kernels $k_s(\cdot, q)\overline{q}u$ are dense in $\mathscr H(s)$.
The isometry property implies that
\[
(0,G)\in R\Longrightarrow G=0.
\]
We can thus introduce a densely defined operator $T$ such that
\[
G=TF\quad\iff\quad (F,G)\in R.
\]
$T$ is a densely defined isometry since $R$ is isometric. As in the case of
complex Hilbert spaces, it extends to an everywhere defined isometry.
We now compute the operator $A$. Let $q\in\mathbb B$ and $u\in
\mathbb H$. We have
\[
A^*(k_s(\cdot, q)\overline{q})u=\left(k_s(\cdot, q)-k_s(\cdot, 0)
\right)u.
\]
Hence, for $f\in\mathscr H(s)$ it holds that:
\[
\begin{split}
\overline{u}q(Af)(q)&=[Af,k_s(\cdot,q)\overline{q}u]_{\mathscr H(s)}\\
&=[f, k_s(\cdot, q)-k_s(\cdot, 0)]_{\mathscr H(s)}\\
&=\overline{u}\left(f(q)-f(0)\right)
\end{split}
\]
and hence
\[
q(Af)(q)=f(q)-f(0).
\]
Similarly we have
\[
B^*(k_s(\cdot, q)\overline{q})u=(\overline{s(q)}-\overline{s(0)})u,
\]
so that we can write for $v\in\mathbb H$:
\[
\begin{split}
\overline{u}q(Bv)(q)&=[Bv,k_s(\cdot, q)u]_{\mathscr H(s)}\\
&=[v,(\overline{s(q)}-\overline{s(0)})u]_{\mathbb H}\\
&=\overline{u}(s(q)-s(0))v
\end{split}
\]
and hence the formula for $B$. To compute $C$ we note that $C^*(\overline{q}v)=
k_s(\cdot, 0)\overline{q}v$ for every $q,v\in\mathbb H$. So,
for $f\in\mathscr H(s)$ we have:
\[
\overline{v}qCf=[Cf,\overline{q}v]_{\mathbb H}=[f,
k_s(\cdot, 0)\overline{q}v]_{\mathscr H(s)}=\overline{v}qf(0),
\]
and so $Cf=f(0)$. Finally, it is clear that $D=s(0)$.
\end{proof}

With these results at hand we turn to the proof of the realization
theorem.\\

{\bf Proof of Theorem \ref{Vaugirard}:}
We note that the pair $(C,A)$ in \eqref{eq:dbh} is closely outer
connected. Let $f\in\mathscr H(s)$, with power series
\[f(p)=\sum_{n=0}^\infty p^nf_n.
\]
We have the formulas
\[
f_n=CA^{n}f,\quad n=0,1,2,\ldots
\]
and
\[
f(p)=C\star(I-pA)^{-\star}f,\quad f\in\mathscr H(s).
\]
Applying these formulas to the function $B1$ we obtain
\[
s(p)-s(0)=pC\star(I-pA)^{-\star}B1.
\]

We now turn to the converse and assume that $s$
is given by a power series which converges in $\mathbb B$ and
for which the coefficients are of the form \eqref{eq:sn}. To prove that
$s$ is a Schur multiplier, we will check the formula
\begin{equation}
\label{eq:Pernety:ligne13}
1-s(p)\overline{s(q)}=U(p)(U(q))^*-pU(p)(U(q))^*\overline{q},\quad p,q\in\mathbb B,
\end{equation}
where
\begin{equation}
\label{eq:U}
U(p)=\sum_{n=0}^\infty p^nCA^n.
\end{equation}

We have
\[
\begin{split}
1-s(p)\overline{s(q)}&=1-(D+\sum_{n=1}^\infty p^nCA^{n-1}B1)
(D+\sum_{m=1}^\infty q^mCA^{m-1}B1)^*\\
&=1-DD^*-\sum_{m=1}^\infty DB^*(A^{m-1})^*C^*\bar q^m-
\sum_{n=1}^\infty p^n CA^{n-1}BD^*+\\
&\hspace{5mm}-\sum_{n,m=1}^\infty p^nCA^{n-1}BB^*(A^{m-1})^*C^*\overline{q}^m\\
&=CC^*+\sum_{m=1}^\infty C(A^{m})^*C^*\bar q^m +
\sum_{n=1}^\infty p^n CA^{n}C^*+\\
&\hspace{5mm}-\sum_{n,m=1}^\infty p^nCA^{n-1}(I-AA^*)
(A^{m-1})^*C^*\overline{q}^m\\
&=U(p)(U(q))^*-pU(p)(U(q))^*\overline{q},
\end{split}
\]
where $U$ is as in \eqref{eq:U} and
we have used the fact that the operator matrix
\eqref{eq:ABCD} is coisometric.\\

It follows from \eqref{eq:Pernety:ligne13} that
\begin{equation}
\label{eq:fondamental}
\sum_{n=0}^\infty
p^n(1-s(p)\overline{s}(q))\overline{q}^n=U(p)U(q)^*,
\end{equation}
is positive definite in $\mathbb B$,
and so  $s$ is a Schur multiplier.\\

Finally we turn to the uniqueness claim. Let
\[
\begin{pmatrix}A_1&B_1\\ C_1&D_1\end{pmatrix}\,\,:\,\,\mathscr
H_1\oplus\mathbb H\longrightarrow \mathscr H_1\oplus\mathbb H
\]
and
\[
\begin{pmatrix}A_2&B_2\\ C_2&D_2\end{pmatrix}\,\,:\,\,\mathscr
H_2\oplus\mathbb H\longrightarrow \mathscr H_2\oplus\mathbb H
\]
be two closely outer-connected coisometric realizations of
$ s$, with state spaces right quaternionic Hilbert spaces
$\mathscr H_1$ and $\mathscr H_2$ respectively. From
\eqref{eq:fondamental} we have
\[
U_1(p)(U_1(q))^*=U_2(p)(U_2(p))^*,
\]
where $U_1$ and $U_2$ are built as in \eqref{eq:U} from the
present realizations. It follows that
\[
C_1A_1^n(A_1^m)^*C_1^*=C_2A_2^n(A_2^m)^*C_2^*,\quad\forall
n,m\in\mathbb N_0.
\]
In view of the presumed outer-connectedness, the relation in
$\mathscr H_1\times \mathscr H_2$ defined by
\[
((A_1^*)^mC_1^*u,(A_2^*)^mC_2^*u),\quad u\in\mathbb H,\quad
m\in\mathbb N_0,
\]
is a densely defined isometric relation with dense range. It is
thus the graph of a unitary map $U$ such that:
\[
U\left((A_1^*)^mC_1^*u\right)=(A_2^*)^mC_2^*u,\quad m\in\mathbb
N_0,\quad{\rm and}\quad u\in\mathbb H.
\]
Setting $m=0$ leads to $UC_1^*=C_2^*$, that is
\begin{equation}
\label{eq1234}
C_1=C_2U.
\end{equation}
With this equality, writing
\[
(UA_1^*)((A_1^*)^mC_1)=A_2^*(A_2^*)^mC_2^*=A_2^*UU^*(A_2^*)^mC_2^*=(A_2^*U)
((A_1^*)^mC_1),
\]
and taking into account that both pairs $(C_1,A_1)$ and
$(C_2,A_2)$ are closely outer-connected, we obtain
$A_1U^*=U^*A_2$, that is
\begin{equation}
\label{eq123} UA_1=A_2U.
\end{equation}
Since clearly $D_1=D_2=s(0)$, it remains only to prove that
$UB_1=B_2$. This follows from the equalities (where we use
\eqref{eq1234} and \eqref{eq123})
\[
s_n=C_1A_1^{n-1}B_1=C_2A_2^{n-1}B_2=C_1A_1^{n-1}U^*B_2, \quad
n=1,2,\ldots
\]
and from the fact that $(C_1,A_1)$ is closely outer connected.
\mbox{}\qed\mbox{}\\

We note that we have followed the arguments in \cite{adrs}
suitably adapted to the present case. In particular the proof of
the uniqueness is adapted from that of \cite[Theorem 2.1.3,  p.
46]{adrs}.\\

The preceding theorem gives a characterization of {\sl all} Schur
multipliers. A simple example is given by the choice
\[
\begin{pmatrix}A&B\\
C&D\end{pmatrix}=\begin{pmatrix}\overline{a}&\sqrt{1-|a|^2}\\
\sqrt{1-|a|^2}&-a\end{pmatrix},
\]
where $a\in\mathbb B$. The corresponding Schur multiplier $s_a(p)$
is
\[
s_a(p)=-a+(1-|a|^2)(1-p\overline{a})^{-\star}=(p-a)\star(1-\overline{a})^{-\star},
\]
and is the counterpart of an elementary Blaschke factor. The
corresponding space $\mathscr H(s_a)$ is finite dimensional and
spanned by the function $(1-p\overline{a})^{-\star}$. We postpone
to a future publication the study of finite dimensional (possibly
indefinite) de Branges Rovnyak spaces in the present setting.\\

We also note that \eqref{eq:Pernety:ligne13} suggests an
equivalent definition of Schur multiplier: The function $s$ is a
Schur multiplier if there is a function $k(p,q)$ positive
definite in $\mathbb B$, and such that
\begin{equation}
\label{sa0}
1-s(p)\overline{s(q)}=k(p,q)-pk(p,q)\overline{q},\quad p,q\in\mathbb B,
\end{equation}
See also \eqref{sa1}, \eqref{sa2}.
These equations, as well as \eqref{sa0}, can be rewritten as sums of
positive definite functions, which induce an isometry. For instance, we can
rewrite
\eqref{sa0} as
\[
1+pk(p,q)\overline{q}=k(p,q)+s(p)\overline{s(q)},\quad p,q\in\mathbb B,
\]
We will not pursue this line of idea here, which is also called
the {\sl lurking isometry method}; see \cite{btv}, and postpone it
to a work where we consider functions of several quaternionic
variables.

\begin{Pn}
The reproducing kernel $k_s$ can be written as
\begin{equation}
\begin{split}
\overline{u}k_s(p,q)v=[(C\star(I-qA)^{-\star})^*v\,
,\,(C\star(I-pA)^{-\star})^*u]_{\mathscr H(s)}\\
&\hspace{-3cm}=\overline {u}\left(\sum_{n=0}^\infty p^n
CA^nA^{*n}C^*\overline{q}^n\right)v.
\end{split}
\end{equation}
\end{Pn}

\begin{proof}
Indeed, both $k_s(p,q)$ and
$U(p)(U(q))^*=\left(\sum_{n=0}^\infty p^n
CA^nA^{*n}C^*\overline{q}^n\right)$ satisfy the equation
\[
X-pX\overline{p}=1-s(p)\overline{s(q)}.
\]
\end{proof}
\section{The Schur algorithm in the quaternionic setting}
\setcounter{equation}{0} \label{schur_h}
Recall that we have denoted by $\mathcal S$ the set of functions
analytic and contractive in the open unit disk $\mathbb D$. We
rewrite the recursion \eqref{recur} in a way which is used in
signal processing, and is conducive to generalization to more
general cases, and in particular to this case. We write
\begin{equation}
\label{classical}
\begin{split}
\begin{pmatrix}1&-s(z)\end{pmatrix}\begin{pmatrix}1&s_0\\
\overline{s_0}&1\end{pmatrix}\begin{pmatrix}z&0\\
0&1\end{pmatrix}&=
\begin{pmatrix}z(1-s(z)\overline{s_0})&-(s(z)-s_0)\end{pmatrix}\\
&=(1-s(z)\overline{s_0})^{-1}\begin{pmatrix}z&-z{s^{(1)}(z)}\end{pmatrix}\\
&=(1-s(z)\overline{s_0})^{-1}z\begin{pmatrix}1&-{s^{(1)}(z)}\end{pmatrix}.
\end{split}
\end{equation}

Let $s(p)=\sum_{p=0}^\infty p^ns_n$ be a Schur multiplier, and
assume that $|s_0|<1$. We now want the counterpart of
\eqref{classical} in the quaternionic setting. We set
\[
s(p)-s(0)=p\sigma^{(1)}(p).
\]
Then, the counterpart of \eqref{classical} in the quaternionic
setting is:
\[
\begin{split}
\begin{pmatrix} I&-M_s\end{pmatrix}\begin{pmatrix}I&M_{s_0}\\
M_{\overline{s_0}}&I\end{pmatrix}\begin{pmatrix}M_p&0\\0&I\end{pmatrix}&=
\begin{pmatrix}(I-M_s M_{\overline{s_0}})M_p&-(M_s-M_{s_0})\end{pmatrix}\\
&=(I-M_s M_{\overline{s_0}})\begin{pmatrix}M_p&-(I-M_s M_{\overline{s_0}})^{-1}
M_pM_{\sigma_1}\end{pmatrix}\\
&=(I-M_s M_{\overline{s_0}})M_p\star\begin{pmatrix}I&-
(I-M_s M_{\overline{s_0}})^{-1}M_{\sigma_1}\end{pmatrix}.
\end{split}
\]

The function $s^{(1)}$ defined by
\[
M_{s^{(1)}}=(I-M_s M_{\overline{s_0}})^{-1}M_{\sigma_1}
\]
will be called the Schur transform of $s$. In view of Theorem
\ref{schwarz_slice} it is a Schur multiplier. When
$|s^{(1)}(0)|<1$, the preceding procedure can be iterated. This
is the quaternionic counterpart of the Schur algorithm.
\section{Concluding remarks}
\setcounter{equation}{0}
Schur analysis has been extended to wide range of other settings.
We mention in particular the definitions of the Schur-Agler
classes for functions analytic in the polydisk or the unit ball.
See \cite{agler-hellinger, ball-trent, BB, btv, MR1909375}. Other
settings include the non-commutative case, \cite{MR2187742,
MR2450144}, Riemann compact surfaces \cite{av3} functions defined
on an homogeneous tree, \cite{MR2481805}, and others. See for
instance \cite{MR93b:47027} for the time-varying case,
\cite{al_acap} for the stochastic case and \cite{am_ieot}
for the case of multi-dimensional systems.\\

We mentioned in the paper a number of problems which will be
considered in future publications. We here list some more:\\

$(1)$ We have presented a coisometric realization of a Schur
multiplier. There are also isometric and unitary realizations,
associated  associated to the reproducing kernel Hilbert spaces
$\mathscr H(s^\sharp)$ and $\mathscr H(D_s)$, where $s^\sharp$
and $D_s$ are defined in \eqref{sharp} and \eqref{eq:DS}
respectively. Their counterpart in the setting of slice
holomorphic functions will be
presented elsewhere.\\

$(2)$ Let $J\in\mathbb H^{n\times n}$ such that $J=J^{*}=J^{-1}$.
A $\mathbb H^{n\times n}$-valued function $\Theta$ defined in an
open subset $\Omega$ of $\mathbb B$ will be called
$J$-contractive if the kernel
\[
K_\Theta(p,q)= \sum_{n=0}^\infty p^n(J-{
\Theta}(p)J({\Theta}(q))^*)\overline{q}^n\\
= (1-\Theta(p)J{(\Theta(q))^*}) \star (1-p\overline{q})^{-\star}
\]
is positive definite. The study of these functions, and their
relations to interpolation problems for Schur multipliers, are an
important topic in Schur analysis. In the classical case, these
functions $\Theta$ are called $J$-contractive, and are
characteristic operator functions; see
\cite{livsic,liv-26,l1,MR48:904,MR20:7221}.
Their structure has been described in \cite{pootapov}.\\

$(3)$ We have mentioned Schur-Agler classes associated to the unit
ball and the poly-disk. Their counterpart in the setting of Fueter
series were presented in \cite{MR2240272,asv-cras}. The study of
Schur multipliers slice hyperholomorphic in several quaternionic
variables should be a topic of interest. Let $p=(p_1,\ldots
,p_N)$ and $q=(q_1,\ldots , q_N)$ vary in $\mathbb B^N$. A
Schur-Agler function is defined as a function $s(p)$ such that
the
\begin{equation}
\label{sa1}
1-s(p)\overline{s(q)}=\sum_{j=1}^N
(k_j(p,q)-p_jk_j(p,q)\overline{q_j}),\quad (\mbox{\rm polydisk
case}),
\end{equation}
or
\begin{equation}
\label{sa2}
1-s(p)\overline{s(q)}=k(p,q)-\sum_{j=1}^N
p_jk(p,q)\overline{q_j},\quad (\mbox{\rm unit ball case}),
\end{equation}
where the $k_j$ and $k$ are positive definite kernels.

\bibliographystyle{plain}
%\bibliography{/users/faculty/math/dany/bib/all}
%\bibliography{all}
\def\cprime{$'$} \def\lfhook#1{\setbox0=\hbox{#1}{\ooalign{\hidewidth
  \lower1.5ex\hbox{'}\hidewidth\crcr\unhbox0}}} \def\cprime{$'$}
  \def\cprime{$'$} \def\cprime{$'$} \def\cprime{$'$} \def\cprime{$'$}

\end{document}